\documentclass[11pt]{article} 
\usepackage{a4,amssymb, amsmath, amsthm}
\usepackage{fullpage}

\newtheorem{thm}{Theorem}[section]
\newtheorem{bigthm}{Theorem}
\newtheorem{cor}[thm]{Corollary}
\newtheorem{lem}[thm]{Lemma}
\newtheorem{pro}[thm]{Proposition}

\theoremstyle{definition}

\theoremstyle{remark}
\newtheorem{rem}[thm]{Remark}

\numberwithin{equation}{section}

\newcommand{\crp}{\overline{\mathbb R}_+}
\newcommand{\rp}{ \mathbb R_+}
\newcommand{\dbar}{d\hspace*{-0.08em}\bar{}\hspace*{0.1em}}
\newcommand{\R}{\mathbb{R}}

\newcommand{\N}{\mathbb{N}}

\newcommand{\Op}{\operatorname{Op}}
\newcommand{\Cal}{\mathcal }
\newcommand{\leg}{\;\dot{\le}\;}
\newcommand{\geg}{\;\dot{\ge}\;}
\newcommand{\eg}{\;\dot {=}\;}
\def\ang#1{\langle {#1} \rangle}

\begin{document}

\title{Heat kernel estimates for pseudodifferential operators, fractional
Laplacians and Dirichlet-to-Neumann operators}
\author{
Heiko Gimperlein\thanks{Department of Mathematical Sciences,  University of Copenhagen, Universitetsparken 5, 2100 Copenhagen \O , 
 Denmark, email: gimperlein@math.ku.dk. Partially supported by the Danish Science
Foundation (FNU) through research grant 10-082866 and the Danish National
Research Foundation through the Centre for Symmetry and Deformation
(DNRF92).
 } \and 
Gerd Grubb\thanks{Department of Mathematical Sciences, University of Copenhagen, Universitetsparken 5, 2100 Copenhagen \O , Denmark, email: grubb@math.ku.dk}}

\maketitle
\begin{abstract}
\noindent The purpose of this article is to establish
 upper and lower estimates for the integral kernel of the
semigroup $\exp(-tP)$ associated to a classical, strongly elliptic
pseudodifferential operator $P$ of positive order on a closed manifold. The Poissonian bounds
generalize those obtained for perturbations of fractional powers of the
Laplacian. In the selfadjoint case, extensions to $t\in{\Bbb C}_+$ are
studied. In particular, our results apply to the Dirichlet-to-Neumann
semigroup.
\end{abstract}

\section*{Introduction}\label{intro}

Let $M$ be a compact $n$-dimensional Riemannian $C^\infty $-manifold and $P$ a
classical, strongly elliptic pseudodifferential operator ($\psi $do) on
$M$ of order $d>0$. We consider upper and lower estimates for the integral
kernel ${\Cal K}_{V}(x,y,t)$ of the generalized heat semigroup
$V(t)=e^{-tP}$. Semigroups generated by such nonlocal operators have been
of recent interest in
different settings.

1) For a Riemannian manifold $\widetilde M$ with boundary $M$, the
Dirichlet-to-Neumann operator is a first-order pseudodifferential operator
on $M$ with principal symbol $|\xi|$.
Arendt and Mazzeo \cite{AM07}, \cite{AM12}, initiated the study of the associated
semigroup and its relation to eigenvalue inequalities, motivating later
studies e.g.\ by Gesztesy and
Mitrea \cite{GM09} and Safarov \cite{S08}.

2) The heat kernel generated by fractional powers of the Laplacian
$\Delta ^{d/2}$ and their perturbations provides another example. Sharp
estimates for $e^{-t\Delta^{d/2}}$, $0<d<2$, can be obtained from those
for $e^{-t\Delta}$ by subordination formulas. For perturbations on bounded
domains in $\Bbb{R}^n$, recent work on estimates includes Chen, Kim and
Song \cite{CKS12} and other works by these authors, and Bogdan et al.\ \cite{BGR10}.

In this article we generalize the Poissonian estimates obtained in the
second case to sectorially elliptic operators $P$ of all positive orders on closed manifolds,
by pseudodifferential methods. In
particular, we allow nonselfadjoint operators and systems. A main
result for such systems $P$ is:
\begin{bigthm} The kernel of the semigroup satisfies
\begin{equation}
\label{intro.thm.1}
|\Cal K_V(x,y,t)|\le C e^{-c_1t}t\,(d(x,y)+t^{1/d})^{-n-d}, \text{ for
  all }x,y\in
M, t\ge 0, \tag{$*$}
\end{equation}
for any $c_1$ smaller than the infimum $\gamma (P)$ of the
real part of the spectrum of $P$.

If $P$ is selfadjoint $\ge 0$, the
estimates extend to the complex $t=e^{i\theta }|t|$ with $|\theta
|<\frac\pi  2$, with uniform estimates
\begin{equation}\label{intro.thm.2}
|\Cal K_{V}(x,y,t)|\le C (\cos \theta  )^{-N}
e^{-\gamma(P)\operatorname{Re}t}\,
\frac{|t|}{(d(x,y)+|t|^{1/d})^{d}}((d(x,y)+|t|^{1/d})^{-n}+1), \tag{$**$}
\end{equation}
where $N=\operatorname{max} \{\tfrac n d , \tfrac{7n}{2} + 4 d + 7 \}$.
\end{bigthm}

Here $d(x,y)$  denotes the distance between $x$ and $y$. If $P$ is a
system, it suffices that $P $ is sectorially elliptic, having
the spectrum of the principal symbol in a sector $\{ \lambda\ne 0 \mid
|\operatorname{arg}\lambda |\le \theta _0\}$ with $\theta _0<\frac
\pi 2$. Extending
\eqref{intro.thm.1}, also derivatives of the kernel, and, if further spectral information is available, a refined
description of the long-time behavior, are obtained in the paper.
Moreover, for the expansion of the kernel in quasi-homogeneous terms
in local coordinates, we show estimates of each term.

For the Dirichlet-to-Neumann operator, as well as for the perturbations of
fractional powers of the Laplacian of orders $0<d<2$, we get not only
upper estimates but also similar lower estimates at small distances.

The estimate \eqref{intro.thm.1} exhibits a large class of operators which satisfy
upper estimates closely related to those studied abstractly e.g.\ in Duong
and Robinson \cite{DR96} and Coulhon and Duong \cite{CD00}. They have implications
for the maximal regularity of the associated evolution equation in $L_p$,
the $L_p$-independence of the spectrum as well as the functional
calculus.

As a simple application of \eqref{intro.thm.2} and H\"o{}lder's
inequality, one can 
obtain ultracontractive estimates
$$
\|e^{-tP}\|_{\Cal L(L_p,L_q)}
\le C (\cos \theta )^{-N}\big\|\frac{|t|}{(d(x,y)+|t|^{1/d})^{d}}((d(x,y)+|t|^{1/d})^{-n}+1)\big\|_{L_{q,y}L_{p',x}},$$
 uniformly for all $t\in{\Bbb C}$ with $\operatorname{Re}t>0$. In the case of operators with
Gaussian heat kernel estimates, a rich spectral theory has been
developed  (see e.g.\ Arendt \cite{A04}, Ouhabaz \cite{O05}).

With the help of comparison principles, our result implies Poissonian
estimates e.g.\ for boundary problems in an open subset $\Omega$ of $M$: If $P$ is the
variational operator associated to a Dirichlet form $a$ with domain $\Cal
D\subset L_2(M)$, we consider the abstract Dirichlet realization
$P_\Omega$ associated to the closure of $a|_{\Cal D \cap C_0(\Omega)}$. In
the case where $a$ is Markovian, one obtains $0\leq \Cal
K_{e^{-tP_\Omega}} \leq \Cal K_{e^{-tP}}$ on $\Omega$. See Grigor'yan and
Hu \cite{GH08} for more refined comparison principles.

\medskip
\noindent {\it Outline.} Section 1 collects some known facts. In Section 2 we
treat semigroups generated by nonselfadjoint $P$ for $t\ge 0$, by
pseudodifferential methods based on \cite{G96}.  Section
3 extends the estimates to complex $t$ for selfadjoint $P$. Section 4
includes lower estimates for perturbations of fractional powers of the
Laplacian and for the Dirichlet-to-Neumann operator.

\section{Preliminaries}

{\it Notation:} $\langle\xi\rangle=\sqrt{\xi ^2+1}$. The indication $\leg$ means ``$\le$ a
constant times'', $\geg$ means  ``$\ge$ a
constant times'', and $\eg$ means that both hold.

\medskip
Let $P$ be a classical $\psi
$do of order $d\in{\Bbb R}_+$, acting in a Hermitian $N$-dimensional
 $C^\infty $
vector bundle $E$ over a closed, compact
Riemannian $n$-dimensional manifold $M$.

We assume that the principal symbol $p^0(x,\xi )$ of $P$ has its spectrum
(for $\xi \ne 0$) in a sector $\{\lambda \ne 0\mid
|\operatorname{arg}\lambda |\le \varphi _0\}$ for some $\varphi _0<\frac
\pi 2$.  (In the notation of the book \cite{G96}, $P-\lambda $ is
parameter-elliptic on the rays in the
complementing sector; according to Seeley \cite{S67}, the latter are ``rays of
minimal growth'' of the resolvent.)
From $P$ one can define the generalized heat operator $V(t)=e^{-tP}$,
$t\ge 0$, a holomorphic semigroup generated by $P$, as explained in
detail e.g.\ in \cite{G96}, Sect.\ 4.2.
The kernel $\Cal
K_{V}(x,y,t)$ ($C^\infty $ for $t>0$)
was analyzed there in its dependence on $t$, but mainly with a view
to sup-norm estimates over all $x,y$, allowing an analysis of the
diagonal behavior, that of $\Cal K_{V}(x,x,t)$. We shall expand the
analysis here to give more information on $\Cal K_{V}(x,y,t)$.
\medskip

For convenience of the reader we recall the definitions of symbol spaces that are used.
For $d\in{\Bbb R}$, the symbol space $S^d_{1,0}({\Bbb R}^n\times{\Bbb
R}^n)$ consists of the $C^\infty $-functions $a(x,\xi )$ ($x,\xi
\in{\Bbb R}^n$) such that for all $\alpha ,\beta \in{\Bbb N}_0^n$,
\begin{equation}\label{1.1}
|D_x^\beta D_\xi ^\alpha a(x,\xi )|\leg \langle\xi\rangle ^{d-|\alpha |};
\end{equation}
it is a Fr\'echet space provided with the seminorms $\sup_{x,\xi
}|\langle\xi\rangle ^{-d+|\alpha |}D_x^\beta D_\xi ^\alpha a|$. The symbols
define operators $A=\operatorname{Op}(a(x,\xi ))$ of order $d$ by
$$
\operatorname{Op}(a(x,\xi ))u= \Cal F ^{-1}
a(x,\xi){\Cal F}u=\int_{{\Bbb R}^n}e^{ix\cdot\xi }a(x,\xi )\hat u(\xi )\, \dbar\xi,
$$
where ${\Cal F} u=\hat u$ denotes the Fourier transform and $\dbar\xi =(2\pi )^{-n}d\xi $. The operator maps from
$\Cal S({\Bbb R}^n)$ to $\Cal S({\Bbb R}^n)$, extending to suitable
spaces of distributions and Sobolev spaces, and obeying various
composition rules.

The space of {\it classical symbols} of order $d$, $S^d({\Bbb R}^n\times{\Bbb
R}^n)$, is the subset of $S^d_{1,0}({\Bbb R}^n\times{\Bbb
R}^n)$ where $a(x,\xi )$ moreover has an asymptotic expansion
$a\sim \sum_{l\in{\Bbb N}_0}a_{d-l}$ in terms $a_{d-l}(x,\xi )$
homogeneous in $\xi $ of degree $d-l$ for $|\xi |\ge 1$, such that
$a'_M=a-\sum_{l<M}a_{d-l}\in S^{d-M}_{1,0}$ for all $M\in {\Bbb
N}_0$. The principal symbol $a_d$ is often denoted $a^0$. (The
homogeneity need only hold for $|\xi |\ge R$, some $R>0$.)

It should be noted that we here use the globally estimated symbols of H\"ormander \cite{H83},
Section 18.1, which have the
advantage that remainders are kept inside the calculus.

Operators on manifolds are defined by use of local coordinates and
rules for change of variables, composition with cut-off functions etc.; we refer to the quoted
works for details.

The book \cite{G96} moreover includes parameter-dependent symbols $a(x,\xi
,\lambda )$ for $\lambda $ in a sector of ${\Bbb C}$, with special
symbol estimates involving the parameter (also operators on manifolds
with boundary are treated there).

Consider a localized situation where the symbol $p(x,\xi)$ of $P$ is defined in a bounded open
subset of ${\Bbb R}^n$ --- we can assume it is extended to
${\Bbb R}^n$, with symbol estimates valid
uniformly in $x$. 
The principal symbol $p^0(x,\xi)$ is an $N\times N$-matrix with spectrum in the
sector $\{\lambda \ne 0\mid |\arg\lambda |\le \varphi _0\}$, when $|\xi |\ge
1$. This holds in particular when $P$ is strongly elliptic, for then
\begin{equation}\label{1.2}
\operatorname{Re}(p^0(x,\xi)v,v)\ge c|\xi|^d|v|^2, \text{ for
}|\xi|\ge 1,\, v\in{\Bbb C}^N,
\text{ with }c>0,
\end{equation}
and hence since
\begin{equation}\label{1.3}
|\operatorname{Im}(p^0v,v)|\le |(p^0v,v)|\le C|\xi|^d|v|^2\le
c^{-1}C\operatorname{Re}(p^0v,v),\text{ for }|\xi|\ge 1,\, v\in{\Bbb C}^N,
\end{equation}
the sectorial ellipticity holds with $\varphi _0=\arctan (c^{-1}C) \in
[0,\frac\pi 2[\,$. When $P$ 
is scalar, the two ellipticity properties are equivalent, but for systems, strong
ellipticity is more restrictive than the mentioned sectorial ellipticity (also
called parabolicity of $\partial_t+P$).

When working in a localized situation, we assume (as we may) that the
sectorial ellipticity holds uniformly for the symbols extended to
$\R^n$. The estimates in the following are valid in particular for operators given
on $\R^n$ with global symbol estimates.

The spectrum $\sigma (P)$ of $P$ lies in a right half-plane and has a finite lower
bound
$\gamma (P)=\inf\{\operatorname{Re}\lambda \mid \lambda \in \sigma (P)\}$.
We can modify $p^0$ for small $\xi$ such that $\sigma (p^0(x,\xi))$
has a positive lower bound for all $(x,\xi )$
 and lies in $
\{\lambda =re^{i\varphi }\mid r>0,\, |\varphi |\le\varphi _0\}$.

The information in the following is taken from \cite{G96}, Section 3.3.

The resolvent $Q_\lambda =(P-\lambda )^{-1}$ exists and is holomorphic in
$\lambda $ on a neighborhood of a set
\begin{equation}\label{1.4}
W_{r_0,\varepsilon }=\{\lambda \in{\Bbb C}\mid |\lambda |\ge r_0,\,
\arg \lambda \in [\varphi _0+\varepsilon ,2 \pi -\varphi _0-\varepsilon
],\}\cup \{ \operatorname{Re}\lambda \le \gamma (P)-\varepsilon \}.
\end{equation}
(with  $\varepsilon >0$).
There
exists a parametrix $Q'_\lambda $
on a neighborhood of a possibly larger set (with $\delta >0,\varepsilon >0$)
$$
V_{\delta ,\varepsilon }=\{\lambda \in{\Bbb C}\mid |\lambda |\ge
\delta \text{ or }
\arg \lambda \in [\varphi _0+\varepsilon ,2 \pi -\varphi _0-\varepsilon
] \}\cup \{\operatorname{Re}\lambda <\inf_{x,\xi }\gamma (p^0(x,\xi ))\};
$$
such that this parametrix coincides with $(P-\lambda )^{-1}$ on
the intersection. Its symbol $q(x,\xi ,\lambda )$ in local coordinates
is holomorphic in $\lambda $ there and has the form
\begin{equation}\label{1.5}
q(x,\xi,\lambda )\sim \sum_{l\ge 0}q_{-d-l}(x,\xi,\lambda )
,\text{ where }
q_{-d}=({p^0(x,\xi)-\lambda })^{-1}.
\end{equation}
Here when $P$ is scalar,
\begin{equation}\label{1.6}
 q_{-d-1}={b_{1,1}(x,\xi
)}q_{-d}^{2},\; \dots,\;
q_{-d-l}=\sum_{k=1}^{2l}{b_{l,k}(x,\xi)}{q_{-d}^{k+1}},\; \dots\; ;
\end{equation}
with symbols $b_{l,k}$  independent of $\lambda $ and homogeneous of
degree   $dk-l$ in $\xi$ for $|\xi|\ge 1$. When $P$ is a system, each
 $q_{-d-l}$ is for $l\ge 1$ a finite sum of terms with the structure
\begin{equation}\label{1.7}
 r(x,\xi ,\lambda )=b_1q_{-d}^{\nu _1}b_2q_{-d}^{\nu _2}\cdots b_Mq_{-d}^{\nu _M}b_{M+1},
 \end{equation}
where the $b_k$ are homogeneous $\psi $do symbols of order $s_k$ independent of
$\lambda $, the $\nu _k$ are positive integers with sum $\in [2, 2l+1]$, and
$s_1+\dots +s_{M+1}-d(\nu _1+\dots+\nu _M)=-d-l$. (Further information
and references in Remark 3.3.7.)
Moreover (cf.\ Theorems 3.3.2 and 3.3.5.), the remainder
$q'_M=q-\sum_{l<M}q_{-d-l}$ satisfies for $\lambda $ on the rays in $W_{r_0,\varepsilon }$: 
\begin{equation}\label{1.8}
|D_{x}^\beta D_{\xi}^\alpha q'_M(x,\xi,\lambda )|\leg \langle\xi\rangle^{d-|\alpha |-M}(1+|\xi|+|\lambda |^{1/d})^{-2d}, \text{ when }M+|\alpha |>d
.\end{equation}

\section{Semigroups generated by sectorially elliptic pseudodifferential
operators}

 As explained in \cite{G96}, Section 4.2, the semigroup $V(t)=e^{-tP }$ can be defined from $P $ by the
 Cauchy integral formula
\begin{equation}\label{2.1}
V(t)=\tfrac i {2\pi }\int_{\Cal C}e^{-t\lambda }(P -\lambda )^{-1}\,d\lambda ,
\end{equation}
where $\Cal C$ is a suitable curve going in the positive direction
around the spectrum of
$P $; it can be taken as the boundary of $W_{r_0,\varepsilon }$ for a
small $\varepsilon $.
 In the local coordinate patch
 the symbol is (for any $M\in{\Bbb N}_0$)
\begin{align}
v(x,\xi,t)&=
v_{-d}+\dots+v_{-d-M+1}+v'_M\sim \sum_{l\ge 0}v_{-d-l}(x,\xi
,t),\text{ where }\nonumber \\
v_{-d-l}&=\tfrac i {2\pi }\int_{\Cal C }e^{-t\lambda }q_{-d-l}(x,\xi
,\lambda )\,d\lambda,\quad v'_{M}=\tfrac i {2\pi }\int_{\Cal C }e^{-t\lambda }q'_{M}\,d\lambda. \label{2.2}
\end{align}

A prominent example is $e^{-t\sqrt\Delta \,}$ where $\Delta $
denotes the (nonnegative) Laplace-Beltrami operator on $M$. This is a  Poisson operator
from $M$ to $M\times\crp$ as defined  in the Boutet de Monvel calculus
(\cite{B71}, cf.\ also \cite{G96}), when $t$ is identified with $x_{n+1}$. When $M$ is replaced by ${\Bbb R}^n$, its kernel is
the well-known Poisson kernel
\begin{equation}\label{2.3}
{\Cal K}(x,y,t)= c_n\frac t{(|x-y|^2+t^2)^{(n+1)/2}}
\end{equation}
for the operator solving the Dirichlet problem
for $\Delta $ on ${\Bbb R}^{n+1}_+$.

Also more general operator families $V(t)=e^{-tP}$ with $P$ of order 1 are sometimes spoken of as Poisson operators (e.g.\ by
Taylor \cite{T81}), and indeed we can show that for $P$ of any order
$d\in\rp$, $V(t)$ identifies with a Poisson operator in the Boutet de
Monvel calculus.
This will be taken up in detail elsewhere. In order to match the conventions for Poisson
symbol-kernels, the indexation in \eqref{2.2} is chosen slightly differently from that in \cite{G96},
Section 4.2, where $v_{-d-l}$ would be denoted $v_{-l}$.
We define $V_{-d-l}(t)$ and $V'_M(t)$ in local coordinates to be the
$\psi $do's with symbol $v_{-d-l}(x,\xi ,t)$ resp.\ $v'_{M}(x,\xi ,t)$.
The
kernel $\Cal K_{V}(x,y,t)$ is in local coordinates expanded according
to the symbol expansion:
\begin{equation}\label{2.4}
\Cal K_V(x,y,t)= \sum_{0\le l<M}\Cal K_{V_{-d-l}}(x,y,t)+\Cal K_{V'_M}(x,y,t).
\end{equation}

The following result follows from \cite{G96}.

\begin{thm}\label{thm2.1}
${\rm 1}^\circ$ In local coordinates,
the kernel terms satisfy for some $c'>0$:
\begin{equation}\label{2.5}
|\Cal K_{V_{-d-l}}(x,y,t)|\leg e^{-c't}\begin{cases} t^{(l-n)/d}\text{ if
}d-l>- n, \\ t\,(
|\log t|+1)\text{ if
}d-l= -n,\\
t\text{ if
}d-l<- n.
\end{cases}
\end{equation}
For a given $c_0>0$ we can modify $p^0$ to satisfy
$\inf_{x,\xi }\gamma (p^0(x,\xi ))\ge c_0$; then $c'$ can be any
number in $\,]0,c_0[\,$.

${\rm 2}^\circ$ Moreover, with the modification in ${\rm 1}^\circ$ used with
$c_0=\gamma (P)$ if $\gamma (P)>0$, the remainder satisfies
\begin{equation}\label{2.6}
|\Cal K_{V'_{M}}(x,y,t)|\leg e^{-c_1t}\begin{cases} t^{(M-n)/d} \text{ if
}d-M>- n, \\ t\,(|
\log t|+1)\text{ if
}d-M= -n,
\\t\text{ if
}d-M<- n,
\end{cases}
\end{equation}
for any $c_1<\gamma (P)$.
In particular,
\begin{equation}\label{2.7}
|\Cal K_{V}(x,y,t)|\leg e^{-c_1t}t^{-n/d}.
\end{equation}
\end{thm}
\begin{proof} The theorem was shown with slightly less precision on the
constants $c',c_1$ in \cite{G96}, Theorems 4.2.2 and 4.2.5. It was there aimed
towards applications where $d$ is integer. The estimates of resolvent
symbols in Section 3.3
are still valid when $d\in\rp$, but the replacement of $P$ by
$P+a$ ($a\in{\Bbb R}$) in the beginning of Section 4.2 on heat operators only gives a
classical $\psi $do when $d$ is integer, so we need another device to
take the value of $\gamma (P)$ into account for general $d\in\rp$. We shall now
explain the needed modifications, with reference to \cite{G96}.

For ${\rm 1}^\circ$, the proof in Theorem 4.2.2 shows the validity of \eqref{2.5}
with a small positive $c'<\inf_{x,\xi }\gamma (p^0(x,\xi ))$. For a given $c_0>0$, the proof goes
through to allow any $c'<c_0$, when $p^0(x,\xi ) $ is modified for $|\xi |\le R$ (for a
possibly large $R$) to satisfy $\inf \gamma (p^0(x,\xi ))\ge
c_0$.

For ${\rm 2}^\circ$, the remainder symbol $q'_M$ is holomorphic on
$W_{r_0,\varepsilon }$; here if $\gamma (P)>0$ we define the terms $q_{-d-l}$ as under
$1^\circ$, with $c_0=\gamma (P)$. For large $M $, $q'_M$ is $\leg \langle\lambda\rangle
^{-2}$. The proof of Th.\ 4.2.2 gives an estimate of $\Cal K_{V'_M}$ by
$e^{-c_{1}t}t\,(1+|\log t|)$, and the proof of Theorem
4.2.5 shows how to remove the logarithm. The estimates of $\Cal K_{V'_M}$ for lower values of $M$ follow by addition of the estimates of finitely many $\Cal K_{V_{-d-l}}$-terms.
\end{proof}

We shall improve this to give information on the dependence on $|x-y|$
also. This will rely on the following result on kernels of
$S^r_{1,0}$-$\psi $do's, found e.g.\ in Taylor \cite{T81}, Lemma XII 3.1, or
\cite{T96}, Proposition VII 2.2.

\begin{pro} \label{pro2.2}
Let $a \in C^\infty(\Bbb R^n \times \Bbb R^n)$
be such that for some $r \in \Bbb R$, some $N \in {\Bbb N}_{0}$ with $N > n+r$,
and all $0 \leq |\alpha|\leq N$,
\begin{equation}\label{2.8}
\operatorname{sup}_{x,\xi } \ang{\xi}^{-r+|\alpha|} |D_\xi^\alpha a(x,
\xi)| \le C_0<
\infty .
\end{equation}
Then the inverse Fourier transform $\Cal K_A(x,y)=\Cal F^{-1}_{\xi \to
z}a(x,\xi )|_{z=x-y}$  is
$O(|x-y|^{-N})$ for $|x-y|\to\infty $, and satisfies
for all $|x-y|>0$:
\begin{equation}\label{2.9}
|\Cal K_A(x,y)|\leg C_0\begin{cases}  |x-y|^{-r-n}\text{ if }r> -n,\\ |\log
|x-y||+1\text{ if }r=-n,\\
1 \text{ if } r<-n.\end{cases}
\end{equation}

In particular, if $a\in S^r_{1,0}({\Bbb R}^n\times{\Bbb R}^n)$ defining the $\psi $do $A$, the estimates hold for its kernel
$\Cal K_A(x,y)$ for all $N>n+r$, each estimate depending only on the
listed symbol seminorms.
\end{pro}

The dependence of the kernel norms on $C_0$ is seen
from an inspection of the proof.

In the scalar case the kernel study can be based on nice
explicit formulas, that we think are worth explaining.
Consider the contribution from one of the terms in \eqref{1.6}. As
integration curve we can here use $C_{\varphi }$ consisting of the two rays
$re^{i\varphi }$ and $re^{-i\varphi }$, $\varphi =\varphi _0+\varepsilon
$. For $t>0$, a replacement of $t\lambda $ by $\varrho $ gives:
\begin{align}\label{2.10}
w_{l,k}(x,\xi,t)&=\tfrac i {2\pi }\int_{C_\varphi  }e^{-t\lambda
}\frac{b_{l,k}(x,\xi)}{(p^0(x,\xi)-\lambda )^{k+1}}\,d\lambda=\tfrac i {2\pi
}\int_{C_\varphi }e^{-\varrho  }\frac{t^{k}b_{l,k}}{(tp^0-\varrho)^{k+1} }\,d\varrho\nonumber \\
&=\tfrac i {2\pi }t^{k}b_{l,k}\int_{C_{\varphi ,R} }\frac{e^{-\varrho
}}{(tp^0-\varrho  )^{k+1}}\,d\varrho =\tfrac 1{k!}t^{k}b_{l,k}e^{-tp^0};
\end{align}
here we have replaced the integration curve by
a closed curve $C_{\varphi ,R}$ connecting the two rays by a circular piece in the
right half-plane  with
radius $R\ge 2t|p^0(x,\xi)|$, and applied the Cauchy integral formula for
derivatives of holomorphic functions. This shows:
\begin{equation}\label{2.11}
v_{-d}=e^{-tp^0},\quad v_{-d-l}(x,\xi,t)=\sum_{k=1}^{2l}\tfrac1{k!}t^k b_{l,k}(x,\xi
)e^{-tp^0(x,\xi)} \text{ for }l\ge 1.
\end{equation}
Then the kernels of the
$V_{-d-l}(t)$ can be estimated by the following observations.

\begin{pro} \label{pro2.3}
Let $p^0(x,\xi )$ be the principal symbol of a
classical scalar strongly elliptic $\psi $do $P$ on ${\Bbb R}^n$ of order
$d\in\rp$, chosen such that $\operatorname{Re}p^0(x,\xi )\ge c_0>0$.

${\rm 1}^\circ$ Let $c'\in[0,c_0[\,$. For any $j\in{\Bbb N}_0$, $(t(p^0(x,\xi )-c'))^je^{-t(p^0(x,\xi )-c')}$ is in
$S^0_{1,0}({\Bbb R}^n\times{\Bbb R}^n)$ uniformly in $t\ge 0$.

${\rm 2}^\circ$
Let
\begin{equation}\label{2.12}w(x,\xi ,t)=\tfrac i {2\pi }\int_{C_\varphi }e^{-t\lambda }\frac{b(x,\xi )}{(p^0(x,\xi )-\lambda )^{k+1}}\,d\lambda,
\end{equation}
where $k\ge 1$ and $b\in S^{dk-l}_{1,0}({\Bbb R}^n\times{\Bbb R}^n)$. Then
\begin{equation}\label{2.13}
w(x,\xi ,t)=\tfrac 1{k!}t^{k}b(x,\xi )e^{-tp^0(x,\xi )}=e^{-c't}t\,w'(x,\xi ,t),
\end{equation}
where 
$w'(x,\xi ,t)
\in S^{d-l}_{1,0}({\Bbb R}^n\times{\Bbb R}^n)$, uniformly  for $t\ge 0$.

Moreover,
$\tilde w(x,z,t)=\Cal F^{-1}_{\xi \to z}w$
satisfies for any $c'\in \,]0,c_0[\,$:
\begin{equation}\label{2.14}
|\tilde w(x,z,t)|\leg e^{-c't}\begin{cases} t\,|z|^{l-d-n}\text{ if }d-l>- n,\\
t\,(|\log |z||+1)\text{ if }d-l= -n ,\\
t\text{ if }d-l< -n .\end{cases}
\end{equation}
It follows that for $l\ge 1$, $\Cal
K_{V_{-d-l}}(x,y,t)=\Cal F^{-1}_{\xi \to z}v_{-d-l}(x,\xi
,t)|_{z=x-y}$ satisfies the estimates
\begin{equation}\label{2.15}
|\Cal K_{V_{-d-l}}(x,y,t)|\leg e^{-c't}\begin{cases} t\,|x-y|^{l-d-n}\text{ if }d-l>- n,\\
t\,(|\log |x-y||+1)\text{ if }d-l= -n ,\\
t\text{ if }d-l< -n .\end{cases}
\end{equation}
Moreover, $\Cal K_{V_{-d-l}}(x,y,t)$ is $O(e^{-c't}t|x-y|^{-N})$ for $|x-y|\to\infty $, any $N$.
\end{pro}
\begin{proof}
${\rm 1}^\circ$. For each fixed $t>0$, $e^{-tp^0(x,\xi )}$ is rapidly decreasing in
$\xi $, hence is in  $S^{-\infty }_{1,0}$. But for our purposes we
need estimates that hold uniformly in $t$ for $t\to 0$. Let
$$
M_{j,k,l}=\sup_{s\ge 0 }s^l\partial_s^k(s^je^{-s}).
$$
Then for all $t\ge 0$, $\xi \in{\Bbb R}^n$,
\begin{align}\label{2.16}
&|(tp^0(x,\xi ))^je^{-tp^0(x,\xi )}|\le M_{j,0,0},\nonumber \\
& |\partial_{\xi _i}\bigl((tp^0 )^je^{-tp^0}
\bigr)|=|\partial_s(s^je^{-s})|_{s=tp^0 }t\partial_{\xi
_i}p^0 |\le M_{j,k,1}|(p^0) ^{-1}\partial_{\xi
_i}p^0 |\leg \ang{\xi }^{-1},\; \dots\nonumber \\
&|\partial_{\xi }^\alpha \bigl((tp^0 )^je^{-tp^0
}\bigr)|\leg \ang{\xi }^{-|\alpha |},\dots
\end{align}
showing the assertion for $c'=0$. \eqref{2.16} holds also if $p^0$ is
replaced by $p^0-c'$ throughout, when $c'\in \,]0,c_0[\,$.

${\rm 2}^\circ$.
The first identity in
\eqref{2.13} was shown in \eqref{2.10}.
We can also write
$$
w(x,\xi ,t)=\tfrac 1{k!}t\,b(p^0-c')^{1-k}(t(p^0-c'))^{k-1}e^{-c't}e^{-t(p^0-c')}=e^{-c't}t\,w'(x,\xi ,t).
$$
Here $b(p^0-c')^{1-k}$ is in $S^{d-l}_{1,0}$, independent of $t$, and by
$1^\circ$, $(t(p^0-c'))^{k-1}e^{-t(p^0-c' )}$ is uniformly in
$S^0_{1,0}$, so it follows that
$w'$ is uniformly in $S^{d-l}_{1,0}$. We can now apply Proposition \ref{pro2.2}
to draw the conclusion \eqref{2.14}.

Since $v_{-d-l}(x,\xi ,t)$ is a sum of such terms when $l\ge 1$, the estimates
\eqref{2.15} follow.

\end{proof}

For systems $P$ we can use systematic estimates from \cite{G96}. We find for
general $P$:

\begin{thm}\label{thm2.4}
${\rm 1}^\circ$ In local coordinates,
$\Cal K_{V_{-d}}$ satisfies for some $c'>0$:
\begin{equation}\label{2.17}
|\Cal K_{V_{-d}}(x,y,t)|\leg e^{-c't} t\, |x-y|^{-d-n}.
\end{equation}
For $l\ge 1$, the kernels $\Cal K_{V_{-d-l}}$ satisfy {\rm \eqref{2.15}}.
For all $l$, the kernels are $O(e^{-c't}t|x-y|^{-N})$ for $|x-y|\to\infty $, any $N$.
If $\gamma (P)>0$,  we  modify $p^0$ to satisfy
$\inf_{x,\xi }\gamma (p^0(x,\xi ))\ge \gamma (P)$, then $c'$ can be any
number in $\,]0,\gamma (P)[\,$.

${\rm 2}^\circ$ Moreover, with $p^0$ chosen as in ${\rm 1}^\circ$,
\begin{equation}\label{2.18}
|\Cal K_{V'_{M}}(x,y,t)|\leg e^{-c_1t}\begin{cases} t\, |x-y|^{M-d-n}\text{ if
}d-M>- n, \\ t\,(|
\log |x-y||+1)\text{ if
}d-M= -n,\\
t\text{ if
}d-M<- n,
\end{cases}\end{equation}
for any $c_1<\gamma (P)$.
In particular,
\begin{equation}\label{2.19}
|\Cal K_{V}(x,y,t)|\leg e^{-c_1t}t\, |x-y|^{-d-n}.\end{equation}
\end{thm}

\begin{proof} ${\rm 1}^\circ$.
When $P$ is scalar, the estimates in \eqref{2.15} for $l\ge 1$ are shown in
Proposition \ref{pro2.3}, when we take $c_0=\gamma (P)$ if $\gamma (P)>0$.
For general systems $P$, the symbols $q_{-d-l}$ are sums of symbols as
in \eqref{1.7}, and we apply \cite{G96}, Lemma 4.2.3. Here (4.2.35) with
$k=-d-l$ shows that
$$
|D_x^\beta D_\xi ^\alpha v_{-d-l}(x,\xi ,t)|\leg \ang\xi ^{d-l-|\alpha |}te^{-c't},
$$
for all $\alpha ,\beta $. Actually, the estimate (4.2.35) has $e^{-ct\ang\xi
^d}$ with a positive $c$ as the last factor, but an inspection of the
proof (the location of integral contours) shows that $e^{-ct\ang\xi
^d}$ can be replaced by $e^{-c't}$,
if $c'<\inf\gamma
(p^0(x,\xi ))$.
This shows that $e^{c't}t^{-1}v_{-d-l}$ is in $S^{d-l}_{1,0}$ uniformly in
$t$, so the estimates of the $\Cal K_{V_{-d-l}}$ follow by use of
Proposition \ref{pro2.2}.

 For $l=0$, we can
argue as follows in the scalar case: For each $j=1,\dots,n$,
$$
\partial_{\xi _j}v_{-d}=\partial_{\xi _j}e^{-tp^0}=-t(\partial_{\xi _j}p^0)e^{-tp^0},$$
where $\partial_{\xi _j}p^0\in S^{d-1}_{1,0}$. Now as in Proposition \ref{pro2.3},
$e^{-c't}\partial_{\xi _j}p^0e^{-t(p^0-c')}$ is in $S^{d-1}_{1,0}$ uniformly in
$t$, and hence $\tilde v_{-d}=\Cal F^{-1}_{\xi \to z}v_{-d}$
satisfies, since $d-1>-n$,
\begin{equation}\label{2.20}
|z_j\tilde v_{-d}|\leg  e^{-c't}t\,|z|^{-d+1-n}.
\end{equation}
Taking the square root of the sum of squares for $j=1,\dots,n$, we find
after division by $|z|$ that
\begin{equation}\label{2.21}
|\tilde v_{-d}|\leg e^{-c't}t\,|z|^{-d-n}.
\end{equation}

In the systems case we note that
\begin{equation}\label{2.22}
\partial_{\xi
_j}q_{-d}=-q_{-d}(\partial_{\xi _j}p^0)q_{-d},
\end{equation}
since
$\partial_{\xi
_j}[(p^0-\lambda )(p^0-\lambda )^{-1}]$
$=0$. Lemma 4.2.3 applies to
this in the same way as above, showing that
$$
|D_x^\beta D_\xi ^\alpha \partial_{\xi _j}v_{-d}(x,\xi ,t)|\leg \ang\xi ^{d-1-|\alpha |}te^{-c't},
$$
so $e^{c't}t\partial_{\xi _j}v_{-d}$ is
is uniformly in $S^{d-1}_{1,0}$. We conclude \eqref{2.20}, from which \eqref{2.21}
follows, implying \eqref{2.17}.

${\rm 2}^\circ$. Here the estimate in \eqref{2.18} has already been shown for large $M$ in
Theorem \ref{thm2.1}.
For lower values of $M$, we can add the
estimates of the entering homogeneous terms $\Cal K_{V_{-d-l}}$ with
$l\ge M$;
the top term gives the weakest estimate. 
\end{proof}

Theorems \ref{thm2.1} and \ref{thm2.4} together lead to Poisson-like kernel estimates:

\begin{thm} \label{thm2.5}
${\rm 1}^\circ$ One has
in local coordinates:
\begin{equation}\label{2.23}
|\Cal K_{V_{-d-l}}(x,y,t)|\leg e^{-c't}\begin{cases} t\,( |x-y|+ t^{1/d})^{l-d-n}\text{ if
}d-l>- n, \\ t\,(|
\log (|x-y|+t^{1/d})|+1)\text{ if
}d-l= -n,\\
t\text{ if
}d-l<- n,
\end{cases}
\end{equation}
for some  $c'>0$. If $\gamma (P)>0$, we  modify $p^0$ to satisfy
$\inf_{x,\xi }\gamma (p^0(x,\xi ))\ge \gamma (P)$; then $c'$ can be any number in $\,]0,\gamma (P)[\,$.

${\rm 2}^\circ$ Moreover, with $p^0$ chosen as in ${\rm 1}^\circ$,
\begin{equation}\label{2.24}
|\Cal K_{V'_{M}}(x,y,t)|\leg e^{-c_1t}\begin{cases} t\,( |x-y|+ t^{1/d})^{M-d-n}\text{ if
}d-M>- n, \\ t\,(|
\log(|x-y|+t^{1/d})|+1)\text{ if
}d-M=  -n,\\
 t\,\text{ if }d-M<  -n,
\end{cases}
\end{equation}
for any $c_1<\gamma (P)$.
In particular,
\begin{align}\label{2.25}
|\Cal K_{V}(x,y,t)|&\leg e^{-c_1t}t\,(|x-y|+t^{1/d})^{-d-n},\nonumber \\
|\Cal K_{V'_1}(x,y,t)|&\leg e^{-c_1t}t\,(|x-y|+t^{1/d})^{1-d-n}.
\end{align}

${\rm 3}^\circ$ For the operators defined on $M$, one has (with $d(x,y)$
denoting the distance between $x$ and $y$)
\begin{equation}\label{2.26}
|\Cal K_{V}(x,y,t)|\leg e^{-c_1t}t\,(d(x,y)+t^{1/d})^{-d-n},
\end{equation}
for any $c_1<\gamma (P)$.
\end{thm}
\begin{proof}  ${\rm 1}^\circ$--${\rm 2}^\circ$. In the region where $|x-y|\ge t^{1/d}$,
\begin{equation*}
|x-y|\le |x-y|+t^{1/d}\le 2|x-y|,
\end{equation*}
in other words, $|x-y|\eg |x-y|+t^{1/d}$. Then the estimates in
Theorem \ref{thm2.4} imply the validity of the above estimates on this region.

In the region where $|x-y|\le t^{1/d}$, we have instead that $t^{1/d}\eg
|x-y|+t^{1/d}$.
Then the estimates in Theorem \ref{thm2.1} imply the above estimates on
that region; for example
$$
t^{-n/d}=t\, (t^{1/d})^{-d-n}\eg t\, (|x-y|+t^{1/d})^{-d-n}
$$
there. For the two regions together, this shows \eqref{2.23}--\eqref{2.25}.

${\rm 3}^\circ$. This follows from the estimates in local coordinates.
\end{proof}

The operator $P$ on $M$ has compact resolvent.
When the eigenvalues with real part equal to $\gamma (P)$
(necessarily finitely many) are semisimple (i.e., the algebraic
multiplicity equals the geometric multiplicity), we can sharpen the
information on the behavior for $t\to\infty $:

\begin{cor} \label{cor2.6}
Assume that all eigenvalues of $P$ with real
part $\gamma(P)$ are semisimple (it holds in particular when $P$ is selfadjoint). Then
\begin{equation}\label{2.27}
|\Cal K_{e^{-tP}}(x,y,t)|\leg e^{-\gamma(P)t}
\frac{t}{(d(x,y)+t^{1/d})^d}\left(
(d(x,y)+t^{1/d})^{-n}+1\right) .
\end{equation}
\end{cor}
\begin{proof}
The spectral projections $\Pi_{j} = \frac{i}{2\pi}\int_{{\Cal C}_j}
(P-\lambda)^{-1} d\lambda$ onto the eigenspaces $X_{j}$ for the
eigenvalues
$\{\lambda _1,\dots,\lambda _k\}$ with real part
$\gamma(P)$ (where $\Cal C_j$ is a small circle around the eigenvalue),
are pseudodifferential operators of order $-\infty$, and their kernels
$\Cal
K_{\Pi _{j}}(x,y)$ are
bounded. If $\varepsilon>0$, the operator $P' = P + \varepsilon
\sum_{j=1}^k\Pi_{j}$ satisfies $\gamma(P')>\gamma(P)$. By
Theorem 2.5 applied to $P'$,
$$
|\Cal K_{e^{-tP'}}(x,y,t)|\leg e^{-\gamma(P)t} t (d(x,y)+t^{1/d})^{-d-n} .
$$
On the other hand, $V(t)=e^{-tP'}+(1- e^{-\varepsilon
t})\sum_{j=1}^ke^{-t\lambda _j}\Pi_{j }$, so
$$
\Cal K_{e^{-tP}}(x,y,t)=\Cal
K_{e^{-tP'}}(x,y,t) +(1- e^{-\varepsilon t})\sum_{j=1}^ke^{-t\lambda _j}
{\Cal K}_{\Pi_{j}}(x,y).
$$
From
$$
1-e^{-\varepsilon t}\leq \min\{1, \varepsilon t\} \leg
\frac{t}{(\operatorname{diam}(M) +t^{1/d})^d} \leq \frac{t}{(d(x,y)
+t^{1/d})^d} ,
$$
we conclude that $(1- e^{-\varepsilon t})|{\Cal K}_{\Pi_j}(x,y)|\leg
\frac{t}{(d(x,y) +t^{1/d})^d}$, and \eqref{2.27} follows since
$|e^{-t\lambda _j}|=e^{-t\gamma (P)}$ for each $j$.
\end{proof}

\begin{rem}\label{rem2.7}
The proof of Corollary \ref{cor2.6} allows to sharpen the estimates in Theorem
\ref{thm2.5} and Theorem \ref{thm2.9} below
also in the general case where the eigenvalues with real part
$\gamma(P)$ are not all semisimple.
Denote by $r$ the dimension of the largest irreducible $P$-invariant
subspace of any eigenspace $X_{j}$ associated to an eigenvalue with real
part $\gamma(P)$. Then in Theorems \ref{thm2.5} and \ref{thm2.9} we may replace the upper bound
$e^{-c't} t (d(x,y)+t^{1/d})^{-d-n-k}$ by
\begin{equation}\label{2.28}e^{-\gamma(P)t}
(1+t^{r-1})\frac{t}{(d(x,y)+t^{1/d})^d}\bigl(
(d(x,y)+t^{1/d})^{-n-k}+1\bigr).
\end{equation}
\end{rem}

It is not hard to extend the estimates to complex $t$ in convenient sectors
around $\rp$. Namely, since $p^0$ has its spectrum in the sector
$\{|\arg \lambda |\le \varphi _0\}$,
$e^{i\theta }P$ has the sectorial ellipticity property when $|\theta |<\theta _0=\frac\pi 2-\varphi _0$.
For each $\theta $ it generates a semigroup $e^{-te^{i\theta }P}$,
and these operator families coincide with the holomorphic extension of
$V(t)$ to the rays $\{re^{i\theta }\}$ in the sector  $V_{\theta _0}=\{t\in{\Bbb C}\mid
|\operatorname{arg}t|<\theta _0\}$. On each ray we have the estimates
in Theorem 2.5, they hold uniformly in closed subsectors of $V_{\theta_0}$. We have hereby obtained:

\begin{thm}\label{thm2.8} 
With $\varphi _0$ defined as in
the beginning of  Section {\rm 1} and $\theta _0=\frac\pi 2-\varphi _0$,
the semigroup
generated by $P$ extends holomorphically to the sector $\{|\arg t|<\theta _0\}$, and the
estimates in Theorem {\rm \ref{thm2.5}} hold in terms of  $|t|$ on any closed sector $\{|\operatorname{arg}t|\le\theta \}$ with
$0<\theta <\theta _0$, taking $c_1<\min_{|\theta '|\le \theta }\gamma (e^{i\theta '}P)$.
\end{thm}

For the case where $P$ is selfadjoint, global estimates on the 
 open sector $\{|\operatorname{arg}t|<\frac\pi 2\}$ will be given in
Section 3.

Also the derivatives
of the kernels can be estimated by use of the symbol estimates in \cite{G96}.

\begin{thm} \label{thm2.9}
${\rm 1}^\circ$ One has
in local coordinates:
\begin{align}
&|D_x^\beta D_{y}^\gamma D_t^j\Cal K_{V_{-d-l}}(x,y,t)|\leg\nonumber \\
& e^{-c't}\begin{cases} t\,( |x-y|+ t^{1/d})^{l-(1+j)d-|\gamma |-n}\text{ if
}(j+1)d+|\gamma | -l>- n, \\ t\,(|
\log (|x-y|+t^{1/d})|+1)\text{ if
}(j+1)d+|\gamma |-l= -n,\\
t\text{ if
}(j+1)d+|\gamma |-l<- n,
\end{cases}\label{2.29}
\end{align}
for some  $c'>0$. If $\gamma (P)>0$, we  modify $p^0$ to satisfy
$\inf_{x,\xi }\gamma (p^0(x,\xi ))\ge \gamma (P)$; then $c'$ can be any number in $\,]0,\gamma (P)[\,$.

${\rm 2}^\circ$ Moreover, with $p^0$ chosen as in ${\rm 1}^\circ$,
\begin{align}
&|D_x^\beta D_{y}^\gamma D_t^j\Cal K_{V'_{M}}(x,y,t)|\leg \nonumber \\
&e^{-c_1t}\begin{cases} t\,( |x-y|+ t^{1/d})^{M-(j+1)d-|\gamma |-n}\text{ if
}(j+1)d+|\gamma |-M>- n, \\ t\,(|
\log(|x-y|+t^{1/d})|+1)\text{ if
}(j+1)d+|\gamma |-M=  -n,\\
 t\,\text{ if }(j+1)d+|\gamma |-M<  -n,
\end{cases}\label{2.30}
\end{align}
for any $c_1<\gamma (P)$.

$3^\circ$ The estimates of derivatives of $\Cal K_V$ hold for the operator defined on $M$ with
$|x-y|$ replaced by $d(x,y)$.
\end{thm}

\begin{proof} As in Theorem \ref{thm2.5}, the estimates are pieced together from
estimates generalizing those in Theorem \ref{thm2.1} resp.\ Theorem \ref{thm2.4} to
include derivatives. We use that
\begin{align*}
&|D_x^\beta D_y^\gamma D_t^j\Cal K_{V_{-d-l}}(x,y,t)|=|D_x^\beta
D_z
^\gamma  D_t^j\tilde v_{-d-l}(x,z ,t)\big|_{z=x-y}|\\
&=|\Cal F^{-1}_{\xi
\to z}(\xi ^\gamma D_x^\beta D_t^j v_{-d-l}(x,\xi ,t))\big|_{z=x-y}|.
\end{align*}

To generalize Theorem \ref{thm2.1} to allow
$x$- and $y$-derivatives we just have to apply the arguments of
\cite{G96}, Theorems 4.2.2 and 4.2.5, to the modified symbols $\xi ^\gamma
D_x^\beta  v_{-d-l}$, to get the  estimates \eqref{2.29}
with  $|x-y|$ replaced by 0. Derivatives with respect to $t$ alone are
explained in Theorem 4.2.5; finally this is combined with $x$- and
$y$-derivatives in a straightforward way.
Similar considerations work for remainders; here we can in fact refer
directly to (4.2.60) for large $M$, and the statements for lower $M$ follow by
addition of the appropriate set of estimates of $\Cal K_{V_{-d-l}}$-terms. This gives the expected generalization of
Theorem \ref{thm2.1}, namely \eqref{2.29}--\eqref{2.30} with $|x-y|$ replaced by 0.

For the generalization of Theorem \ref{thm2.4} we note that estimates
$$
|\xi ^\gamma D_x^\beta D_\xi ^\alpha D_t^jv_{-d-l}(x,\xi ,t)|\leg
\ang\xi ^{(j+1)d+|\gamma |-|\alpha |-l}t e^{-c't}
$$
for $|\alpha  |+l>0$, all $\beta $, $j$, follow
from \cite{G96}, Lemma 4.2.3
(see the remarks around (4.2.40) for how to include
$t$-derivatives, as done also in Theorem 4.2.5). Thus
$e^{c't}t^{-1}\xi ^\gamma D_x^\beta D_t^jv_{-d-l}$ is in
$S^{(j+1)d+|\gamma |-l}_{1,0}$ uniformly in $t$, and it follows by Proposition
\ref{pro2.2} that
\begin{equation*}
|D_z ^\gamma D_x^\beta D_t^j\tilde v_{-d-l}(x,z ,t)|\leg
 e^{-c't}\begin{cases} t|z|^{-(j+1)d-|\gamma |+l-n},\text{ if
 }(j+1)d+|\gamma | -l>- n, \\
t\,(|
\log (|z|+t^{1/d})|+1)\text{ if
}(j+1)d+|\gamma |-l= -n,\\
t\text{ if
}(j+1)d+|\gamma |-l<- n.\end{cases}
\end{equation*}
This implies
estimates as in \eqref{2.29} with $|x-y|+t^{1/d}$ replaced by $|x-y|$.
The conclusion is immediate for $l\ge 1$, and for $l=0$, we use the
estimates of $D_{\xi _j}v$ as in the proof of Theorem \ref{thm2.4}.
Again for remainder estimates, we can appeal to (4.2.60) for
large $M$.

The proof is now completed as in Theorem \ref{2.5}.
\end{proof}

\begin{rem}
As an example of a non-selfadjoint strongly elliptic case of interest
to which the results apply, let us
mention the first-order operator $P=\Delta ^{\frac12}+L$, where $L$ is a first-order
differential operator with real coefficients (in the situation on $\R^n$, $L=b(x)\cdot \nabla
+c(x)$, $b$ and $c$ real smooth with all derivatives bounded). Since the principal symbol $b(x)\cdot i\xi $ of $L$ is purely imaginary,
$\operatorname{Re}p^0(x,\xi )=|\xi |$, so $p^0(x,\xi )$ for $\xi \ne 0$ ranges in a
sector $\{\lambda \ne 0\mid |\arg\lambda |\le \varphi _0\}$, $\varphi _0<\frac\pi
2$. This case is treated by other methods in Xie and Zhang \cite{XZ12}. 

\end{rem}

\section{Estimates in the complex plane for selfadjoint operators}

In this section we shall derive some uniform kernel estimates for the
extension of the semigroup into the region ${\Bbb C}_+=\{t
\in{\Bbb C}\mid \operatorname{Re}t >0\}$, when $P$ is
selfadjoint. We assume for simplicity that $P$ is $\ge 0$.
As noted in Theorem \ref{2.8}, $V(t)$ exists for all $t\in {\Bbb
C}_+$, and the estimates worked out in Section 2 hold uniformly on
closed subsectors  $\{ t\in {\Bbb C}_+\mid |\arg t|\le \theta  \}$,
$0\le\theta  <\frac\pi 2$. For the analysis of the behavior for
$\theta \to \frac\pi 2$, additional efforts are needed. We shall 
rely on a theorem of Agmon \cite{A65}:

\begin{pro}\label{pro5.1}
 Let $\Omega $ be an open set having the cone
property. Let $T$ be a bounded operator in $L_2(\Omega )$ such that
the ranges of $T$ and $T^*$ are contained in $H^m(\Omega )$ for an
$m>n$ ($m$ can be noninteger if $\Omega ={\Bbb R}^n$ or is suitably
smooth). Then $T$ is an integral operator with a continuous and
bounded kernel $\Cal K_T(x,y)$ on $\Omega \times\Omega $ satisfying
$$
|\Cal K_T(x,y)|\le C(\|T\|_{0,m}+\|T^*\|_{0,m})^{n/m}\|T\|_{0,0}^{1-n/m},\label{5.1}
$$
with a constant depending only on $\Omega $ and $m$.
\end{pro}

Here $\|T\|_{a,b}$ stands for the norm of $T$ as an operator from
$H^a(\Omega )$ to $H^b(\Omega )$. The theorem holds in particular when
$\Omega $ is replaced by our manifold $M$ or by $\R^n$.

We note first that there is an easy estimate of the kernel in terms of
$t$ alone, 
that can be obtained essentially by functional analysis.

\begin{thm}\label{thm5.2} 
The full kernel satisfies for $t=e^{i\theta }|t|\in{\Bbb C}_+$:
$$
|{\Cal K}_{V}(x,y,t)|\leg 1+(\cos\theta)^{-n/d}|t|^{-n/d}  .
$$
\end{thm}
\begin{proof}
Note that $\|V(t)\|_{0,m} \eg \|(1+P^{m/d})V(t)\|_{0,0}$ because
$P\geq0$ is strongly elliptic and of order $d$. By the standard theory
of analytic semigroups (see e.g.\ Lunardi \cite{L95}, Proposition 2.1.1) or functional calculus, $((\operatorname{Re} t)P)^M V(t)$ is uniformly bounded for $M \in \N_0$. This extends to $M \geq 0$ by interpolation. Hence we obtain for $m\geq 0$:
$$\|V(t)\|_{0,m} \eg \|(1+P^{m/d})V(t)\|_{0,0} \leg
1+(\operatorname{Re} t)^{-m/d}\ .$$
The estimates likewise hold for $V(t)^*=V(\bar t)$. 
Proposition \ref{pro5.1} therefore yields for any $m>n$
$$|{\Cal K}_{V}(x,y,t)|\leg (\|V(t)\|_{0,m}+\|V(\bar t)\|_{0,m})
^{n/m}\|V(t)\|_{0,0}^{1-n/m} \leg  1+(\operatorname{Re} t)^{-n/d}
\eg 1+(\cos\theta\,|t|)^{-n/d}\ , $$
as was to be shown.\qed
\end{proof}

Extensions to $\Bbb C_+$ of the other estimates in Section 2 
are more costly in negative powers of  $\cos\theta $. We first
consider the
homogeneous terms in the symbol of $Q_\lambda  $,
showing how the estimates of
symbols  like \eqref{1.5}--\eqref{1.7} depend  on $\arg
\lambda $, when $\lambda $ is close to the spectrum of $P$.

As in \cite{G96}, we denote $|\lambda |^{1/d}=\mu $, and write $\ang{(\xi ,\mu )}=(1+|\xi
|^2+\mu ^2)^{1/2}$ for short as $\ang{\xi ,\mu }$;
it is $\eg (1+|\xi |+|\lambda |^{1/d})$.

\begin{pro} \label{pro3.1}
Let $p^0(x,\xi )$ be symmetric with lower bound $\ge c\ang{\xi }^d$
, and let
$\lambda \in{\Bbb C}$ with $\arg \lambda =\varphi $, $0<|\varphi |\le \frac\pi 2$. Then
\begin{align}
|q_{-d}(x,\xi ,\lambda )|&=|(p^0(x,\xi )-\lambda )^{-1}|\leg |\sin \varphi |^{-1}\ang{\xi ,\mu }^{-d},\nonumber \\
|D_{x}^\beta D_{\xi}^\alpha q_{-d}(x,\xi,\lambda )|&\leg |\sin \varphi
|^{-1-|\alpha |-|\beta |}\ang\xi ^{d-|\alpha |}\ang{\xi ,\mu }^{-2d}
,\text{ when }|\alpha |+|\beta |>0. \label{3.1}
\end{align}
For all $l,\alpha ,\beta $ with $l>0$,
\begin{equation}\label{3.2}
|D_{x}^\beta D_{\xi}^\alpha q_{-d-l}(x,\xi,\lambda )|\leg |\sin \varphi |^{-2l-1-|\alpha |-|\beta |}\ang\xi ^{d-l-|\alpha |}\ang{\xi ,\mu }^{-2d}.
\end{equation}
\end{pro}
\begin{proof}
We have for  $\lambda =e^{i\varphi }|\lambda |$ with $0<|\varphi |\le
\frac\pi 2$, and  $v\in{\Bbb C}^N$:
\begin{align*}
|(p^0v,v)-\lambda |v|^2|&\ge |\operatorname{Im}((p^0v,v)-|\lambda
|e^{i\varphi }|v|^2)|=|\lambda ||\sin\varphi ||v|^2,\\
|(p^0v,v)-\lambda |v|^2|&=|e^{-i\varphi }(p^0v,v)-|\lambda ||v|^2|\ge
|\operatorname{Im}e^{-i\varphi }(p^0v,v)|\\&=|\sin\varphi |(p^0v,v)\ge
|\sin\varphi |\,c\ang \xi ^d|v|^2,
\end{align*}
from which follows
$$
|(p^0-\lambda )v||v|\ge |((p^0-\lambda )v,v)|\geg |\sin\varphi |(|\lambda |+\ang\xi ^d)|v|^2.
$$
 This implies that $|(p^0-\lambda )^{-1}|\leg |\sin\varphi
 |^{-1}(\ang\xi ^d+|\lambda |)^{-1}\eg |\sin\varphi
 |^{-1}\ang{\xi ,\mu }^{-d}$, showing \eqref{3.1}.

The other estimates follow as in \cite{G96} from the structure of the
terms in the parametrix, using \eqref{3.1}: $q_{-d-l}$ is for $l\ge 1$ a
finite sum of terms,
where $\nu _1+\dots+\nu _M\ge 2$ takes
values up to $2l+1$,
$$
r(x,\xi ,\lambda )=b_1q_{-d}^{\nu _1}b_2q_{-d}^{\nu _2}\cdots b_Mq_{-d}^{\nu _M}b_{M+1},
$$
cf.\  \eqref{1.7}. Each $q_{-d}$ contributes to the estimates with a
factor $|\sin\varphi |^{-1}$,
and there are up to $2l+1$ such factors; this shows \eqref{3.2} for $\alpha
=\beta =0$. Each differentiation may hit a factor $q_{-d}$ giving an
extra $|\sin\varphi |^{-1}$ in view of \eqref{2.22}; this leads to the
estimates \eqref{3.2} by the Leibniz formula.
\end{proof}

The symbol terms $v_{-d-l}(x,\xi ,t)$ are defined from the $q_{-d-l}(x,\xi ,\lambda )$ as in \eqref{2.2}.
When
$t\in {\Bbb C}_+$, with argument $\arg t=\theta
\in \,]-\frac\pi 2,\frac\pi 2[$,   we must assure that
$\operatorname{Re}(\lambda t)\to \infty $ when $|\lambda |\to\infty $ on the integral curve $\Cal C$. This holds
if $\lambda $ runs on a contour formed of the rays $\lambda =re^{\pm
  i\varphi _0}$, where $\varphi _0=\frac12(\frac\pi 2-|\theta |)$,
connected near 0 by a circle of radius $\varepsilon '<\inf\gamma (p^0(x,\xi )) $
passing  to the
right of 0:
\begin{equation}
\Cal C=\{re^{i\varphi _0}\mid \infty >r>\varepsilon '\}\cup
\{\varepsilon 'e^{i\varphi }\mid \varphi _0>\varphi >-\varphi _0
\}\cup\{re^{-i\varphi _0}\mid \varepsilon '<r<\infty \},\quad \varphi _0=\tfrac12(\tfrac\pi 2-|\theta |).\label{3.1a}
\end{equation}
Here $\inf _{\lambda \in\Cal C}\operatorname{Re}(\lambda )
=c_1>0$.

\begin{rem} \label{rem3.4}
Note that $|\theta |=\frac\pi 2-2\varphi _0$  belongs to $[0,\frac\pi
2[\,$ if and only if $\varphi _0$ belongs to
$\,]0,\frac\pi 4]$. On this interval, $\sin\varphi _0\eg \sin(2\varphi
_0)=\cos\theta $, so they can be used interchangeably in our estimates.
\end{rem}

We shall need the following generalization of \cite{G96} Lemma 4.2.3.

\begin{lem}\label{lem3.9} Let $t=e^{i\theta }|t|$, and choose $\varphi
 _0 $ and $\Cal C$ as in {\rm \eqref{3.1a}}.

Let $M\in\Bbb N$, let $\sigma_1,\dots,\sigma_{M}$ be nonnegative integers with
\begin{equation}
 \sigma=\sigma_1+\cdots+\sigma_{M}\ge 1,\label{3.29}
\end{equation}
and let $f(x,\xi,\lambda)$ be a (matrix-formed) symbol
of the form
\begin{equation}
f(x,\xi,\lambda)=f_1(p^0-\lambda)^{-\sigma_1}f_2(p^0-\lambda)^{-\sigma_2}
\cdots(p^0-\lambda)^{-\sigma_M}f_{M+1},\label{3.30}
\end{equation}
 where the $f_j(x,\xi)$ are $\psi $do symbols of order $s_j\in\Bbb
R$, homogeneous for $|\xi|\ge1$.
 Denote $s_1+\cdots+s_{M+1}=s$, then the order of $f$ is $k=s-\sigma d$. Let
$F_\lambda=\Op(f(x,\xi,\lambda))$ on $\Bbb
R^n$, and let $E(t)$ be the operator family defined from $F_\lambda$
for $\operatorname{Re}t>0$ by
 \begin{equation}
E(t)=\tfrac i{2\pi}\int_{\Cal C}e^{-t\lambda}F_\lambda
\,d\lambda.\label{3.31}
 \end{equation}
 Then $E(t)=\Op(e(x,t,\xi))$, where the symbol
 \begin{equation}
e(x,t,\xi)=\tfrac
i{2\pi}\int_{\Cal C}e^{-t\lambda}f(x,\xi,\lambda)\,d\lambda\label{3.32}
 \end{equation}
satisfies:
 \begin{eqnarray}
&\text{\rm (i)} &\quad
e(x,s^{-d}t,s\xi)=s^{d+k}e(x,t,\xi)\text{ for }
|\xi|\ge 1,\;s\ge 1,\nonumber\\
& \text{\rm (ii)} &\quad
|D_x^\beta D_\xi^\alpha e(x,t,\xi)|\leg (\sin \varphi_0)^{-\sigma
-|\alpha|-|\beta|}
\ang\xi ^{d+k-|\alpha|}e^{-c\operatorname{Re}t \ang \xi ^d}.\label{3.33}
 \end{eqnarray}
The kernel of $E(t)$ satisfies for $d+k>- n$
\begin{equation}|\Cal K_E(x,y,t)| \leg (\sin \varphi_0)^{-\sigma-(d+k+n)/d}
e^{-c'\operatorname{Re}t} |t|^{-(d+k+n)/d} ,\label{3.34}\end{equation}
with $c'>0$.

If $\sigma \geq 2$,
\begin{equation}|D_x^\beta D_\xi^\alpha e(x,t,\xi)| \leg  (\sin
\varphi_0)^{-\sigma-|\alpha|-|\beta|} |t| \ang\xi
^{2d+k-|\alpha|}e^{-c'\operatorname{Re}t \ang \xi ^d} \ . \label{3.35}\end{equation}
In this case, the kernel satisfies for $d+k\leq - n$
\begin{equation}|\Cal K_E(x,y,t)| \leg (\sin \varphi_0)^{-\sigma}
e^{-c'\operatorname{Re}t} \begin{cases}
|t|\,(
|\log \operatorname{Re}t |+1)\text{ if
}d+k= -n,\\
 |t|\text{ if
}d+k<- n.
\end{cases} \label{3.36}\end{equation}
\end{lem}
\begin{proof}
As in \cite{G96}, Lemma 4.2.3, we can pass the operator definition through
the integral.
To estimate $e$, we first consider $|\xi| \leq 1$. We use the residue
theorem and that $p^0$ is selfadjoint to obtain
\begin{equation*}|e(t,x,\xi)| = |\tfrac i{2\pi}\int_{\Cal C}e^{-t\lambda}f(x,\xi,\lambda)\,d\lambda| \leg (1+|t|^{\sigma-1})
e^{-c\operatorname{Re}t}\ .\end{equation*}
Here, $c = \gamma(p^0(x,\xi))$.

For $|\xi| \geq 1$, we replace $\Cal C$ by a closed,
homogeneous curve $\Cal C_{c,C}$ around the spectrum of
$p^0(x,\xi)$. $\Cal C_{c,C}$ coincides with
$\Cal C$ on a annulus of inner radius $c |\xi|^{d}$ and
outer radius $C |\xi|^{d}$ and is closed by the segments of the boundary of
this annulus which lie to the right of $\Cal C$. Then by
homogeneity,
$$|e(t,x,\xi)|= |\tfrac i{2\pi}\int_{\Cal C_{c,C}}e^{-t\lambda}f(x,\xi,\lambda)\,d\lambda| \leg (\sin
\varphi_0)^{-\sigma} \ang \xi ^{d} \ang \xi ^{k} e^{-\tfrac c 2
\operatorname{Re}t |\xi|^d}\ .$$
Combining the two estimates, we conclude \eqref{3.33} for $\alpha = \beta = 0$.
The derivatives $D_x^\beta D_\xi^\alpha e(x,t,\xi)$ are sums of terms of a similar form,
with $k$ replaced by $k-|\alpha|$ and $\sigma$ replaced by numbers
$\le \sigma+|\alpha|+|\beta|$.

To show \eqref{3.34} for $d+k>- n$, we estimate $\Cal K_E$ by comparing $e$ with its
homogeneous extension $e^h$:
$$\Cal K_E(x,y,t) = \int_{\Bbb R^n} e^{i (x -y)\cdot \xi} e^h(x,t,\xi)\,
\dbar\xi + \int_{|\xi|\leq 1} e^{i (x -y)\cdot \xi} (e-e^h)\, \dbar\xi \ .$$
Using \eqref{3.33} and a homogeneous variant,
\begin{align*}
|\Cal K_E(x,y,t)| &\leg (\sin \varphi_0)^{-\sigma} e^{-c_1
\operatorname{Re}t} \int_{\Bbb R^n} e^{-c_2 \operatorname{Re}t
|\xi|^{d}} |\xi|^{d+k}\, \dbar\xi \\
& \quad + (\sin \varphi_0)^{-\sigma} e^{-c_1
\operatorname{Re}t} \int_{|\xi|\leq 1} e^{-c_2 \operatorname{Re}t
|\xi| ^{d}}(\ang \xi ^{d+k} + |\xi|^{d+k})\, \dbar\xi \ .
\end{align*}
The first integral is $\eg (\operatorname{Re}t)^{-(d+k+n)/d} \eg (\sin \varphi_0)^{-(d+k+n)/d} |t|^{-(d+k+n)/d}$, while the second remains
bounded as $|t|\to 0$.

Now consider the case where $\sigma \geq 2$. As $|f|\leg \ang\lambda
^{-2}$ away from $\Bbb R_+$, the integral converges uniformly in $t\geq
0$. We may deform $\Cal C$ to a closed curve in the left
half-plane, where $f$ is holomorphic, to conclude $e(x,0,\xi) = 0$. Also,
using $(-\lambda)(p_d-\lambda)^{-1} = 1-p_d(p_d-\lambda)^{-1}$,
\begin{equation*}
\partial_t e (x,t,\xi)= \tfrac
i{2\pi}\int_{\Cal C}e^{-t\lambda}(-\lambda)
f(x,\xi,\lambda)\,d\lambda
\end{equation*}
can be expressed in terms of $e$ and a second term of the same form, with
one of the $s_j$ replaced by $s_j+d$. By \eqref{3.33}
$$
|\partial_t e(x,t,\xi)| \leg (\sin \varphi_0)^{-\sigma} \ang\xi
^{2d+k}e^{-c\operatorname{Re}t \ang \xi ^d}
$$
and hence, since the value at $t=0$ is 0,
$$|e(x,t,\xi)| \leg  (\sin \varphi_0)^{-\sigma} |t| \ang\xi
^{2d+k}e^{-c\operatorname{Re}t \ang \xi ^d} \ .$$
This shows \eqref{3.35} for $\alpha = \beta = 0$. The proof for $D_x^\beta
D_\xi^\alpha e(x,t,\xi)$ is analogous. The estimate \eqref{3.36} is obtained similarly to \eqref{3.34}, using \eqref{3.35} instead of \eqref{3.33}.
\end{proof}

This leads to the estimates of homogeneous terms:

\begin{thm}\label{thm3.10} 
Let $t=e^{i\theta }|t|\in{\Bbb C}_+$.
In local coordinates, the homogeneous terms  in the kernel of $V(t)$ satisfy for some $c'>0$:
\begin{align}
|\Cal K_{V_{-d}}(x,y,t)|&\leg (\cos\theta )^{-n/d}e^{-c'\operatorname{Re}t}|t|^{-n/d},\\
|\Cal K_{V_{-d-l}}(x,y,t)|&\leg (\cos\theta )^{-2l-1}
e^{-c'\operatorname{Re}t}\begin{cases} (\cos\theta )^{(l-n)/d}|t|^{(l-n)/d}\text{ if
}d-l>- n, \\  |t|\,(
|\log \operatorname{Re}t |+1)\text{ if
}d-l= -n,\\
 |t|\text{ if
}d-l<- n.
\end{cases}\label{3.37}
\end{align}
\end{thm}
\begin{proof} We choose
$\varphi_0 $ and the curve $\Cal C$ as in \eqref{3.1a}, recalling that 
  $\cos\theta \eg \sin\varphi_0 $. 
 For $l \geq 1$, the assertion follows from Lemma \ref{lem3.9},
\eqref{3.34} resp.~\eqref{3.36}, using that in the terms of
$q_{-d-l}$, $(d+k+n)/d=(n-l)/d$ and $\sigma \le 2l+1$.

For $l=0$, we explicitly compute 
\begin{align*}
|\Cal K_{V_{-d}}&(x,y,t)| = (2 \pi)^{-n}|\int_{\Bbb R^n} e^{i (x -y)\cdot \xi}
e^{-tp^0(x,\xi)}\, \dbar\xi |\\
& \leg e^{-c_1 \operatorname{Re}t} \big(\int_{\Bbb R^n} e^{-c_2
\operatorname{Re}t |p^0_h(x,\xi)|}\, \dbar\xi + \int_{|\xi|\leq 1} (e^{-c_2 \operatorname{Re}t |p^0(x,\xi)|}-e^{-c_2
\operatorname{Re}t |p^0_h(x,\xi)|})\, \dbar\xi\big)\\
& \leg e^{-c_1 \operatorname{Re}t} \big(\int_{\Bbb R^n} e^{-\operatorname{Re}t |\xi|^d}\, \dbar\xi
+ 1\big) \leg e^{-c' \operatorname{Re}t} (\operatorname{Re}t)^{-n/d} \ .
\end{align*}
The assertion follows, since $\operatorname{Re}t = |t|\cos \theta $.

\end{proof}

Moreover, estimates in terms of $|t|$ and powers of $|x-y|$ 
are obtained. For $a\in{\Bbb R}$ we denote by $[a]$ the largest
integer $\le a$.

\begin{thm}\label{thm3.11}
${\rm 1}^\circ$ In local coordinates,
$\Cal K_{V_{-d}}$ satisfies for some $c'>0$:
\begin{equation}
|\Cal K_{V_{-d}}(x,y,t)|\leg (\cos\theta)^{-[d-1+n]-3} e^{-c'\operatorname{Re}t}
|t|\, |x-y|^{-d-n}.\label{3.38}
\end{equation}
For $l\ge 1$, the kernels $\Cal K_{V_{-d-l}}$ satisfy
\begin{align}
& |\Cal K_{V_{-d-l}}(x,y,t)|\leg \nonumber\\
& (\cos\theta)^{-2l-1}
e^{-c'\operatorname{Re}t}\begin{cases} (\cos\theta)^{-[d-l+n]-1}|t|\,|x-y|^{l-d-n}\text{ if }d-l>- n,\\
(\cos\theta)^{-1} |t|\,(|\log |x-y||+1)\text{ if }d-l= -n ,\\
 |t|\text{ if }d-l< -n .\end{cases}
\label{3.39}
\end{align}
${\rm 2}^\circ$ Moreover, 
\begin{align}
& |\Cal K_{V_{-d-l}}(x,y,t)|\leg 
e^{-c'\operatorname{Re}t}\begin{cases} (\cos\theta)^{-N_l}|t|\,(|x-y|+|t|^{1/d})^{l-d-n}\text{ if }d-l>- n,\\(\cos\theta)^{-N_l}
|t|\,(|\log (|x-y|+|t|^{1/d})|+1)\text{ if }d-l= -n ,\\
(\cos\theta)^{-N_l} 
 |t|\text{ if }d-l< -n ,\end{cases}
\nonumber\\
&\text{where }N_l=\begin{cases}\max\{n/d,[d-1+n]+3\}\text{ if }l=0,\\
\max\{2l+1+(n-l)/d,2l+2+[d-1+n]\}\text{ if }l>0,d-l>-n,\\
2l+2\text{ if }d-l=-n,\\
2l+1\text{ if }d-l<-n.\\
\end{cases}\label{3.39a}
\end{align}

\end{thm}

\begin{proof}
${\rm 1}^\circ$. For $l\geq 1$, we obtain from Lemma \ref{lem3.9}, \eqref{3.35}, that
$$|D_\xi^\alpha v_{-d-l}(x,t,\xi)| \leg (\cos\theta)^{-2l-1-|\alpha|} \ang
\xi ^{d-l-|\alpha|} |t| e^{-c'\operatorname{Re}t}\ .$$
Here we apply Proposition \ref{pro2.2} with $r=d-l$; then we need $|\alpha
|\le N$ where $N\in{\Bbb N}_0$, $N> d-l+n$. If $d-l+n<0$ we take $N=0$, and
if $d-l+n\ge 0$, we take $N=[d-l+n]+1$. This shows \eqref{3.39}.

For $l=0$, $v_{-d}(x,t,\xi) = e^{-tp^0(x,\xi)}$, and we pass via
$\partial_{\xi _j}v$ to show the estimate, as in Theorem \ref{thm2.4}. 
Here $\partial_{\xi
  _j}v$ enters by application of Lemma \ref{lem3.9} to $f(x,\xi ,\lambda )=\partial_{\xi
  _j}q_{-d}=-q_{-d}(\partial_{\xi _j}p^0)q_{-d}$; it has $\sigma =2$,
$k=-d-1$. Then \eqref{3.35} gives
that
$$|D_\xi^\alpha \partial_{\xi _j}v_{-d}(x,t,\xi)| \leg (\cos\theta)^{-2-|\alpha|} \ang
\xi ^{d-1-|\alpha|} |t| e^{-c'\operatorname{Re}t}\ ,$$ 
so an application of  Proposition \ref{pro2.2} with $N=[d-1+n]+1$ gives that
$$
|z_j\tilde v_{-d}|\leg (\cos\theta)^{-2-[d-1+n]-1} |t|
e^{-c'\operatorname{Re}t}|z|^{-d+1-n}.$$ 
Using this for all $j=1,\dots, n$, we find \eqref{3.38}.

$2^\circ$. We here combine the preceding estimates with those in
Theorem \ref{thm3.10} in the same way as in the proof of Theorem \ref{thm2.5}. 
\end{proof}

Estimates of remainders $V_M'$ are more
difficult to
work out, since they depend on the interplay between the exact
resolvent $Q_\lambda $ and the homogeneous symbol terms, and they will
be more
costly in powers of $(\cos\theta )^{-1}$, the larger $M$ is taken. We
shall here go directly to remainder {\it kernel} estimates.

Let us define the $M$-th resolvent remainder
operator $$
Q'_M=Q_\lambda -\sum_{l<M}Q_{-d-l},
$$ where each $Q_{-d-l}$ is an
operator on the manifold $M$ constructed from the symbols $q_{-d-l}$
in local coordinates. For each $\lambda $, $Q'_M$ is a $\psi $do of
order $-d-M$, but we do not know on beforehand how it is estimated in terms
of $\lambda $, although we have such information on the terms
$Q_{-d-l}$.
Let us write
\begin{align}
Q'_M&=Q'_M(P-\lambda )Q_\lambda = R_MQ_\lambda ,\text{ where
}\nonumber \\
 R_M&=Q'_M(P-\lambda )=1-{\sum}_{l<M}Q_{-d-l}(P-\lambda )\label{3.3}
\end{align}
is a $\psi $do of order $-M$ constructed from known symbols. 
The idea is now that functional analysis gives us a certain control over operator norms of
$Q_\lambda $, whereas $\psi $do calculus will allow us to estimate
operator norms of $ R_M$, and then Agmon's result Proposition
\ref{pro5.1} will lead to a kernel estimate of the composed operator.
For $\lambda $ with argument $\varphi $ satisfying $0<|\varphi |\le
\frac\pi 2$, 
\begin{align}
\|(P-\lambda )u\| \|u\|&\ge |((P-\lambda )u,u)|\ge
|\operatorname{Im}\lambda |\|u\|^2=|\sin\varphi |\,|\lambda
|\,\|u\|^2;\text{ hence }\nonumber\\
\|Q_\lambda \|_{0,0}&\le |\sin\varphi |^{-1}\,|\lambda |^{-1}.\label{3.2b}
\end{align}

We are aiming for an estimate of $\Cal K_{V'_M}$ by $c|t|$, and we
know from Theorem \ref{thm2.1} based on \cite{G96}, Thm.\ 4.2.5, that
to avoid logarithmic factors it is better to use the resolvent formula
\begin{equation}
Q_\lambda = -\lambda ^{-1}+\lambda ^{-1}Q_\lambda P,\label{3.3a}
\end{equation}
inserted in the integral \eqref{2.1} derived with respect to $t$.

First some details on how to handle the possible zero eigenspace of
$P$.
Similarly to
Corollary \ref{cor2.6}, it will be convenient to write $P=P^\varepsilon
-\varepsilon \Pi _0$, where $\Pi _0$ is the orthogonal projection onto
the zero eigenspace of $P$, and $\varepsilon >0$ is chosen $\le $ the
lowest positive eigenvalue, whereby $P^\varepsilon =P+\varepsilon
\Pi _0$ is $\ge \varepsilon $. Here $\Pi _0$ is the $\psi $do of order 0 with kernel ${\sum}_{j=1}^\nu
\varphi _j(x)\varphi _j(y)^*$, for an orthonormal basis $\varphi
_1,\dots,\varphi _\nu $ of the zero eigenspace.
Then $$
V(t)=V^\varepsilon (t)+(1-e^{-\varepsilon t})\Pi _0,\text{ where }V^\varepsilon (t)=e^{-tP^\varepsilon };
$$
and it is the latter operator that needs investigation. $V^\varepsilon (t)$ is defined from the resolvent
$Q^\varepsilon _\lambda =Q_\lambda -(\varepsilon -\lambda )^{-1}\Pi
_0=(P^\varepsilon -\lambda )^{-1}$ by
\begin{equation}
V^\varepsilon (t)=\tfrac i {2\pi }\int_{\Cal C}e^{-t\lambda }Q^\varepsilon _\lambda \,d\lambda .\label{3.19}
\end{equation}
For this integral, the contour can be chosen as in \eqref{3.1a} with $\varepsilon '<\varepsilon $.

For simplicity of notation we drop the $\varepsilon $-index in the
next calculations, and return to include the contribution from $\Pi _0$ in the
final formulations.

An application of \eqref{3.3a} gives
\begin{align}
V(t)&=\tfrac i {2\pi }\int_{\Cal C}e^{-t\lambda }(-\lambda
^{-1}+\lambda ^{-1}Q _\lambda P)\,d\lambda =\tfrac i {2\pi }\int_{\Cal
C}e^{-t\lambda }\lambda ^{-1}Q _\lambda P\,d\lambda,\nonumber\\
\partial_tV(t)&=-\tfrac i {2\pi }\int_{\Cal C}e^{-t\lambda }Q _\lambda P\,d\lambda.
\label{3.20}
\end{align}
Thus $\partial_t\Cal K_{V'_M}$ is the kernel of the integral of the  $M$-th
remainder $-(Q_\lambda P)'_M$ of $-Q _\lambda P$. We know that $\Cal K_{V'_M}$
vanishes at $t=0$, and want to show boundedness of the last integral
applied to the kernel of the $M$-th remainder.
Here (cf.\ also \eqref{3.3})
$$
Q_\lambda P=({\sum_{l<M}}Q_{-d-l}+Q'_M)P=
{\sum}_{l<M}Q_{-d-l}P+  R_MQ _\lambda P={\sum}_{l<M}Q_{-d-l}P+  R_MPQ _\lambda .
$$
Since $R_MPQ_\lambda $ is already of order $-M$, the
 $M$-th remainder of $Q_\lambda P$ is
\begin{align}
(Q_\lambda P)'_M=\widetilde R_M+R_MPQ_\lambda , \text{ where
}\widetilde R_M=({\sum}_{l<M}Q_{-d-l}P)'_M.\label{3.3b}
\end{align}

In preparation for the study of the symbols of $R_M$ and $\widetilde
R_M$,  we prove a lemma on composition 
formulas from the $\psi $do theory. In the general first two rules it is important that the
component to the right is $\lambda
$-independent, to keep the introduction of factors $|\sin\varphi|^{-1}$ as low as possible.

\begin{lem}\label{lem3.2}
Let
$b(x,\xi ) \in
S^{d_2}_{1,0}({\Bbb R}^n\times{\Bbb R}^n)$, and let $a(x,\xi ,\lambda ) \in S^{d_1}_{1,0}({\Bbb R}^n\times{\Bbb R}^n)$
with respect to $(x,\xi )$, with $\lambda $ as in Proposition {\rm
\ref{pro3.1}}, such that for some $d'\ge 0$, $N \in \Bbb R$, one has for all $\alpha ,\beta \in{\Bbb N}_0^n$,
\begin{equation}\label{3.5}
|D_x^\beta D_\xi ^\alpha a(x,\xi ,\lambda )|\leg
|\sin\varphi|^{-N-|\alpha|-|\beta|}\ang\xi^{d' + d_1-|\alpha|}
\ang{\xi ,\mu }^{-d'} \ .
\end{equation}

${\rm 1}^\circ$ There exists $c (x,\xi ,\lambda ) \in S^{d_1+d_2}_{1,0}({\Bbb R}^n\times{\Bbb R}^n)$
such that $\operatorname{Op}(a) \operatorname{Op}(b) =
\operatorname{Op}(c)$, and for every $M \in{\Bbb N}_0$, 
\begin{equation}\label{3.6}
c(x,\xi ,\lambda ) = {\sum}_{|\alpha|<M} \tfrac{1}{\alpha!}\ D_\xi^\alpha
a(x,\xi ,\lambda )
\partial_x^\alpha b(x,\xi) + c_M(a,b) \ ,
\end{equation}
where
\begin{equation}\label{3.7}
|D_x^\beta D_\xi ^\alpha c_M(a,b)| \leg
|\sin\varphi|^{-N-M-|\alpha|-|\beta|}\ang\xi^{d'+d_1+d_2-M-|\alpha|}
\ang{\xi ,\mu }^{-d'}.
\end{equation}

${\rm 2}^\circ$ If {\rm \eqref{3.5}} for $\alpha =\beta =0$ is replaced by
\begin{equation}\label{3.8}
|a(x,\xi ,\lambda )|\leg
|\sin\varphi|^{-N}\ang\xi^{d + d_1}
\ang{\xi ,\mu }^{-d} \ .
\end{equation}
for some $0\le d\le d'$, then {\rm \eqref{3.6}} holds with
{\rm \eqref{3.7}} valid for $M\ge 1$ and the
estimates of $c_0$  replaced by
\begin{equation}\label{3.9}
|D_x^\beta D_\xi ^\alpha c_0(a,b)| \leg
|\sin\varphi|^{-N-|\alpha|-|\beta|}\ang\xi^{d+d_1+d_2-|\alpha|}
\ang{\xi ,\mu }^{-d}.
\end{equation}
${\rm 3}^\circ$ For $\gamma \in{\Bbb N}_0^n$, $D^\gamma \Op(a)=\Op(a^\gamma)
$, where
\begin{equation}\label{3.10}
|D_x^\beta D_\xi ^\alpha a^\gamma (x,\xi ,\lambda )|\leg
{\sum}_{k\le |\gamma |}|\sin\varphi|^{-N-k-|\alpha|-|\beta|}\ang\xi^{d' + d_1+|\gamma |-k-|\alpha|}
\ang{\xi ,\mu }^{-d'} \ .
\end{equation}
\end{lem}
\begin{proof} ${\rm 1}^\circ$. Let $\chi (x,\xi )$ denote a $C^\infty
$-function that is 1 for $|x|^2+|\xi |^2\le 1$ and vanishes for
$|x|^2+|\xi |^2\ge 2$, then we can replace the given symbols by
their products with $\chi (\varepsilon x,\varepsilon \xi )$, which
makes all integrals calculated below convergent.
It is known in the theory (by the technique of oscillatory integrals, cf.\ \cite{H83},
Sect.\ 7.8), that the resulting symbols converge to the given symbols
for $\varepsilon \to 0$
in all the seminorms that are involved. The modified symbols will again
be denoted $a$, $b$.
We can also assume that $b$ has compact support in $x$
(in a set containing the $x$ for which we need the formula). Then
$
\hat b(\eta ,\xi )=\Cal F_{x\to \eta }b(x,\xi  )
$
 satisfies
\begin{equation}\label{3.11}
|D_\xi ^\alpha \hat b(\eta ,\xi )|\leg \ang\eta ^{-N'}\ang\xi ^{d_2-|\alpha |},
\end{equation}
for all $\alpha , N'$.
It follows from the $\psi $do
defining  formula
that $\operatorname{Op}(a)\Op(b)=\Op (c)$, where
\begin{align}
c(x,\xi ,\lambda )&=\int_{{\Bbb R}^{4n}}a(x,\eta ,\lambda
)b(y,\xi )e^{i(x-y)\cdot\eta }e^{i(y-z)\cdot\xi }\,dz\dbar\xi dy\dbar\eta \nonumber \\
&=\int_{{\Bbb R}^n}a(x,\xi +\eta ,\lambda)\hat b(\eta ,\xi
)e^{ix\cdot\eta }\,\dbar\eta .\label{3.12}
\end{align}
If $M>0$, we insert the Taylor expansion of $a$ in $\xi $ up to order $M$,
\begin{align*}
a(x,\xi +\eta ,\lambda )&= {\sum}_{|\alpha |<M}\tfrac1{\alpha !}\eta
^\alpha \partial_{\xi }^\alpha a(x,\xi ,\lambda )\\
& \qquad +{\sum}_{|\alpha
|=M}\tfrac M{\alpha !}\eta ^\alpha \int_0^1(1-h)^{M-1}\partial_\xi ^\alpha a(x,\xi +h\eta ,\lambda )\,dh,
\end{align*}
obtaining that $c=c_{<M}+c_M$, where
\begin{align*}
c_{<M}&=(2\pi )^{-n}\int_{{\Bbb R}^n}{\sum}_{|\alpha |<M}\tfrac1{\alpha
!}\partial_\xi ^\alpha a(x,\xi ,\lambda )\eta ^\alpha \hat b(\eta ,\xi
)e^{ix\cdot\eta }\,\dbar\eta\\
&={\sum}_{|\alpha|<M} \tfrac{1}{\alpha!}\partial_\xi^\alpha
a(x,\xi ,\lambda )
D_x^\alpha b(x,\xi)={\sum}_{|\alpha|<M} \tfrac{1}{\alpha!}D_\xi^\alpha
a(x,\xi ,\lambda )
\partial_x^\alpha b(x,\xi), \\
c_M&=(2\pi )^{-n}\int_{{\Bbb R}^n}{\sum}_{|\alpha |=M}\tfrac M{\alpha
!}\int_0^1(1-h)^{M-1}\partial_\xi ^\alpha a(x,\xi +h\eta ,\lambda )\,dh \,\eta ^\alpha \hat b(\eta ,\xi
)e^{ix\cdot\eta }\,\dbar\eta.
\end{align*}
The sum over $|\alpha |<M$ equals the sum in \eqref{3.6}. For the
last integral we use that
\begin{align*}
|\partial_\xi ^\alpha a(x,\xi +h\eta ,\lambda )|&\leg |\sin\varphi
|^{-N-M}\ang{\xi +h\eta }^{d'+d_1-M}\ang{\xi +h\eta ,\mu }^{-d'}\\ &\leg
|\sin\varphi
|^{-N-M}\ang{\xi }^{d'+d_1-M}\ang{\xi ,\mu }^{-d'}\ang\eta ^{|d'+d_1-M|+d'},
\end{align*}
by the Peetre inequality. Taking this together with the estimates
\eqref{3.11} of $\hat b$ (with a large $N'$), we can conclude
that
\begin{equation*}
|c_M|\leg
|\sin\varphi
|^{-N-M}\ang{\xi }^{d'+d_1+d_2-M}\ang{\xi ,\mu }^{-d'}.
\end{equation*}

For $M=0$, we apply such considerations directly to $c_0=c(x,\xi ,\lambda
)$ in \eqref{3.12}:
$$
|c_0|\leg |\sin\varphi |^{-N}\int \ang{\xi +\eta }^{d'+d_1}\ang{\xi
+\eta ,\mu }^{-d'}\ang\eta ^{-N'}\ang\xi ^{d_2}\,\dbar\eta \leg
|\sin\varphi |^{-N} \ang{\xi  }^{d'+d_1+d_2}\ang{\xi
 ,\mu }^{-d'}.
$$

Derivatives of $c_M$ in $x$ and $\xi $ are treated in a similar way.

In the case ${\rm 2}^\circ$ the proof goes through in a similar way, except
that $d'$ is replaced by $d$ in expressions containing
undifferentiated factors $a$.

In ${\rm 3}^\circ$,  the $\lambda $-independent factor is to the left, and
\eqref{3.12} holds with integrand \linebreak$(\xi +\eta )^\gamma \Cal F_{z\to\eta }a(z,\xi
,\lambda )e^{ix\cdot\eta }$. The Taylor expansion of $(\xi +\eta )^\gamma $ is
a finite binomial expansion ${\sum}_{\kappa \le \gamma }\tbinom
\gamma  \kappa \xi ^{\gamma -\kappa }\eta ^\kappa $ and leads to  a finite
composition formula where the estimates \eqref{3.10} of the terms can be read off directly.
\end{proof}

The composed symbol $c=c_0(a,b)$ is also denoted $a\circ b$ (used in
\cite{G96}) or $a\# b$.

For the analysis of $ R_M$, we denote $P-\lambda =\widetilde P$, with the
parameter-dependent symbol $\tilde p(x,\xi ,\lambda )=p(x,\xi
)-\lambda $ in local coordinates; here for any $M\in {\Bbb N}_0$,
\begin{align}
p&=\sum_{k<M}p_{d-k}+p'_M,\quad \tilde p=\sum_{k<M}\tilde
p_{d-k}+\tilde p'_M,\text{ with}\nonumber \\
\tilde p_d&=p-\lambda ,\quad \tilde p_{d-k}=p_{d-k}\text{ for
}k>0,\quad \tilde p'_M=p'_M\text{ for }M>0.\label{3.13}
\end{align}
$p_d$ is also denoted $p^0$.
The $p_{d-k}$ are homogeneous in $|\xi |$ of degree $d-k$ for $|\xi
|\ge 1$, and $p'_M\in S^{d-M}_{1,0}({\Bbb R}^n\times{\Bbb R}^n)$.

\begin{pro} \label{pro3.3}
Let $M\ge 1$. 
The symbol $ r_M(x,\xi ,\lambda )$ of
$ R_M$ (cf.\ {\rm \eqref{3.3}}) satisfies in local coordinates:
\begin{equation}\label{3.14}
|D_x^\beta D_\xi ^\alpha  r_M(x,\xi ,\lambda )|\leg |\sin\varphi
|^{-2M-|\alpha |-|\beta |}\ang\xi ^{d-M-|\alpha |}\ang{\xi ,\mu
}^{-d}.
\end{equation}
We also have that $R_M=R^{(1)}_M+ R^{(2)}_M$ with symbols 
\begin{equation} r_M= r^{(1)}_M+r_M^{(2)},\quad r_M^{(2)}=q_{-d }p'_M,\label{3.14a}
\end{equation}
estimated by:
\begin{align}
|D_x^\beta D_\xi ^\alpha  r^{(1)}_M(x,\xi ,\lambda )|&\leg
|\sin\varphi |^{-2M-|\alpha |-|\beta |}\ang\xi ^{2d-M-|\alpha
  |}\ang{\xi ,\mu }^{-2d},\nonumber\\
|D_x^\beta D_\xi ^\alpha  r^{(2)}_M(x,\xi ,\lambda )|&\leg |\sin\varphi |^{-1-|\alpha |-|\beta |}\ang\xi ^{d-M-|\alpha |}\ang{\xi ,\mu }^{-d}.\label{3.14b}
\end{align}
Moreover, $\widetilde R_M =-R_M$.
\end{pro}

\begin{proof}
We have that
$$
 r_M=1-{\sum}_{k<M}{\sum}_{l<M}q_{-d-l}\circ
\tilde p_{d-k}- {\sum}_{l<M}q_{-d-l}\circ
\tilde p'_{M}.
$$

The terms in the parametrix symbol ${\sum}_{l\geq
0} q_{-d-l}$ are constructed as  solutions to the successive equations
for $m\in{\Bbb N}_0$:
\begin{equation}
{\sum}_{|\alpha| + k + l = m}
\tfrac{1}{\alpha!}\ D_\xi^\alpha q_{-d-l} \partial_x^\alpha \tilde p_{d-k} =\begin{cases}
 1 \text{ for }m=0,\\
0\text{ for }m =1,2,\dots,
\end{cases} \label{3.15}
\end{equation}
cf.\ e.g.\ Seeley \cite{S67}, (1). We use the truncated composition formula in Lemma \ref{lem3.2}
to compute the
symbol ${r}_{  M}\ $ of ${R}_{  M}\
$ with expansions in up to  ${M}$ homogeneous terms:
\begin{align*}
{r}_{  M} &=
1-{\sum}_{k<M}{\sum}_{l<M}
\Bigl\{{\sum}_{k+l+|\alpha| < {M}} \tfrac{1}{\alpha!}\ D_\xi^\alpha
q_{-d-l} \partial_x^\alpha  \tilde p_{d-k} + c_{{M}-k-l}(q_{-d-l},
\tilde p_{d-k})\Bigr\} \\
 & \qquad - c_0\Bigl({\sum}_{l<M}q_{-d-l},
\tilde p'_{M}\Bigr) \ .
\end{align*}
By \eqref{3.15},
$$
{\sum}_{k<M}{\sum}_{l<M}
{\sum}_{k+l+|\alpha| <M} \tfrac{1}{\alpha!}\ D_\xi^\alpha q_{-d-l}
\partial_x^\alpha \tilde p_{d-k} = 1.
$$
Thus $ r_M$ consists of the following terms:
\begin{align}
{r}_{  M}  =
-{\sum}_{k<M}{\sum}_{l<M}
c_{{M}-k-l}(q_{-d-l}, \tilde p_{d-k}) - c_0\Bigl({\sum}_{l<M}q_{-d-l}, \tilde p'_{M}\Bigr) \ .\label{3.15a}
\end{align}

Using the estimates \eqref{3.1}  and \eqref{3.2} together with
$|D_x^\beta D_\xi^\alpha \tilde p_{d-k}(x,\xi)| \leg \ang\xi^{d-k-|\alpha|}$,
we obtain from Lemma \ref{lem3.2} with $d'=2d$, $d_1 = -d-l$ and $d_2
= d-k$ that for $l \geq 1$ in the  sum over $k,l$:
\begin{align}
|D_x^\beta &D_\xi^\alpha c_{{M}-k-l}(q_{-d-l}, \tilde p_{d-k})|\nonumber \\
&\leg
|\sin\varphi|^{-{M}+k+l-1-2l-|\alpha|-|\beta|}\ang\xi^{2d-d-l+d-k-({M}-k-l)-|\alpha|}
\ang{\xi ,\mu }^{-2d}\nonumber \\
& \leq 
|\sin\varphi|^{- 2M -|\alpha|-|\beta|}\ang\xi^{2d-{M}-|\alpha|}\ang{\xi ,\mu }^{-2d},
\label{3.16}
\end{align}
since $k-l\ge -M+1$.

For $l=0$ we find in view of Lemma \ref{lem3.2} ${\rm 2}^\circ$, since $k<M$,
\begin{align}
|D_x^\beta &D_\xi^\alpha c_{{M}-k}(q_{-d}, \tilde p_{d-k})|\nonumber \\
 &\leg
|\sin\varphi|^{-1-{(M-k)}-|\alpha|-|\beta|}\ang\xi^{2d-d+d-k-{(M-k)}-|\alpha|}\ang{\xi
,\mu }^{-2d}\nonumber  \\
&\leg
|\sin\varphi|^{-{M}-1-|\alpha|-|\beta|}\ang\xi^{2d-{M}-|\alpha|}\ang{\xi ,\mu }^{-2d}
\ .
\label{3.17}
\end{align}

In the last term,
\begin{align}
|D_x^\beta &D_\xi^\alpha  c_0\Bigl({\sum}_{1\le l<M}q_{-d-l},
\tilde p'_{M}\Bigr)| \nonumber \\
 &\leg 
{\sum}_{1\le l<M}|\sin\varphi|^{-2l-1
-|\alpha|-|\beta|}\ang\xi^{2d-d-l+(d-M)
-|\alpha|}\ang{\xi ,\mu }^{-2d}\nonumber  \\
& \leg |\sin\varphi|^{-2  M
-|\alpha|-|\beta|}\ang\xi^{2d-
M-|\alpha|}\ang{\xi ,\mu }^{-2d} 
\label{3.18},
\end{align}
whereas
$$
c_0(q_{-d}, p'_M)=q_{-d}p'_M+c_1(q_{-d}, p'_M),
$$
with $c_1(q_{-d}, p'_M)$ estimated as in \eqref{3.18} and 
\begin{equation}
|D_x^\beta D_\xi^\alpha (q_{-d}p'_M)|\leg  |\sin\varphi|^{-1-|\alpha|-|\beta|}\ang\xi^{d-
M-|\alpha|}\ang{\xi ,\mu }^{-d}.\label{3.18a}
\end{equation}

An addition of the contributions (using that $\ang\xi /\ang{\xi ,\mu
}\le 1$) gives \eqref{3.14}. We also have the representation
\eqref{3.14a}, where all the contributions to
$r^{(1)}_M$ have $O(\ang{\xi ,\mu }^{-2d})$, so that it satisfies
\eqref{3.14b},
and $r^{(2)}_M$ is estimated in \eqref{3.18a}.

For the analysis of $\widetilde R_M$ we  have by Lemma \ref{lem3.2}:
\begin{align}
{\sum}_{l<M}q_{-d-l}\circ
 p&={\sum}_{k<M}{\sum}_{l<M}q_{-d-l}\circ
 p_{d-k}+ {\sum}_{l<M}q_{-d-l}\circ
 p'_{M}\nonumber\\
&={\sum}_{k<M}{\sum}_{l<M}
\Bigl\{{\sum}_{k+l+|\alpha| < {M}} \tfrac{1}{\alpha!} D_\xi^\alpha
q_{-d-l} \partial_x^\alpha   p_{d-k} + c_{{M}-k-l}(q_{-d-l},
 p_{d-k})\Bigr\}\nonumber\\
&\quad + c_0\Bigl({\sum}_{l<M}q_{-d-l},
 p'_{M}\Bigr),\nonumber
\end{align} 
such that the $M$-th remainder has symbol
\begin{align}
\tilde r_M
&={\sum}_{k<M}{\sum}_{l<M}
 c_{{M}-k-l}(q_{-d-l},
 p_{d-k})+ c_0\bigl({\sum}_{l<M}q_{-d-l},
 p'_{M}\bigr)\nonumber\\
&={\sum}_{k<M}{\sum}_{l<M}
c_{{M}-k-l}(q_{-d-l},
 p_{d-k})+ c_0\bigl({\sum}_{1\le l<M}q_{-d-l},
 p'_{M}\bigr)+q_{-d}p'_M+c_1(q_{-d}, p'_M).\nonumber
\end{align}

Here we can observe that all the $p$-factors can be replaced by the
corresponding $\tilde p$-factors, for they are the same when the index
is $\ne d$, and $p_d$ enters only in
differentiated form since $l<M$ (and 
 $\tilde p_d=p_d-\lambda $ and $p_d$ have the same derivatives).
Then in view of the formula \eqref{3.15a} for $r_M$, we have indeed
$\tilde r_M=-r_M$ and $\widetilde R_M=-R_M$.

\end{proof}

Summing up, we now have (cf.\ \eqref{3.3b} ff.), since $\widetilde R_M=-R_M$,
$$
(Q_\lambda P)'_M=- R_M+R_MPQ_\lambda =-R^{(1)}_M-R^{(2)}_M+R^{(1)}_MPQ_\lambda +R^{(2)}_MPQ_\lambda .
$$
These terms will enter in different ways in the integral defining $V'_M(t)$.
Some further symbol estimates will be needed in the following:

\begin{pro}\label{pro3.6} For $\gamma \in {\Bbb N}_0^n$, $k\in {\Bbb N}_0$, the symbols
  of $R^{(1)}_M\ang D^k $, $D^\gamma R^{(1)}_M $, $R^{(2)}_MP\ang D^k $ and
  $D^\gamma R^{(2)}_MP $ satisfy:
\begin{align}
|D_x^\beta D_\xi ^\alpha (r^{(1)}_M\circ \ang {\xi }^k )|&\leg  |\sin\varphi |^{-2M-|\alpha |-|\beta |}\ang\xi ^{2d-M-|\alpha
  |+k}\ang{\xi ,\mu }^{-2d},\nonumber\\
|D_x^\beta D_\xi ^\alpha (\xi ^\gamma \circ r^{(1)}_M)|&\leg  |\sin\varphi |^{-2M-|\alpha |-|\beta |-|\gamma |}\ang\xi ^{2d-M-|\alpha
  |+|\gamma |}\ang{\xi ,\mu }^{-2d},\nonumber\\
|D_x^\beta D_\xi ^\alpha (r^{(2)}_M\circ p\circ \ang {\xi }^k )|&\leg  |\sin\varphi |^{-1-|\alpha |-|\beta |}\ang\xi ^{2d-M-|\alpha
  |+k}\ang{\xi ,\mu }^{-d},\nonumber\\
|D_x^\beta D_\xi ^\alpha (\xi ^\gamma \circ r^{(2)}_M\circ p)|&\leg  |\sin\varphi |^{-1-|\alpha |-|\beta |-|\gamma |}\ang\xi ^{2d-M-|\alpha
  |+|\gamma |}\ang{\xi ,\mu }^{-d}.\label{3.21}
\end{align}

\end{pro}

\begin{proof}
The composition with $\ang D^k =\operatorname{Op}(\ang\xi ^k)$ to the right just corresponds to
multiplying the symbol by $\ang{\xi }^k $, so the first line in
\eqref{3.21} results directly from \eqref{3.14b}. For the second line
we use the composition rule in Lemma \ref{lem3.2} $3^\circ$.
For the third line, we can for $k =0$ use the composition rule in
Lemma \ref{lem3.2} $1^\circ$, which gives the result in view of
\eqref{3.14b}. The third line with $k \ne 0$ follows simply by
multiplication by $\ang{\xi }^k$, and the fourth line follows by
another application of Lemma \ref{lem3.2} $3^\circ$.
\end{proof}

To apply Agmon's estimate Proposition \ref{pro5.1} to obtain kernel estimates, we need an estimate of
$L_2$-bounds in terms of symbol seminorms. Many variants are known,
and we use the following, found in Marschall \cite{M87}, Theorem 2.1.

\begin{pro}\label{pro3.5}
Let $a \in S^0_{1,0}(\Bbb R^n \times \Bbb R^n)$ be such that for some
$C_0>0$, some $N \in {\Bbb N}_{0}$ with $N>\frac{n}{2}$, and all $\alpha, \beta
\in {\Bbb N}_{0}^n$ with $0 \leq |\alpha|
\leq N$, $0 \leq |\beta|\leq 1$,
\begin{equation}
\operatorname{sup}_{x,\xi } \ang{\xi}^{-|\alpha|} |D_x^\beta D_\xi^\alpha
a(x, \xi)| \leq C_0 <
\infty .\label{3.21}
\end{equation}
Then the associated operator $A = \operatorname{Op}(a)$ is bounded on
$L_2(\Bbb R^n)$, and $\|A\|_{0,0} \leg C_0$.
\end{pro}

The dependence of the operator norm on $C_0$ follows from an inspection
of the proof.

\begin{thm}\label{thm3.5a}
The kernels of the operators $-R^{(1)}_M$, $R^{(1)}_MPQ_\lambda $ and
$R^{(2)}_MPQ_\lambda$  are estimated by
\begin{align}
|{\Cal K}_{R^{(1)}_M}(x,y,\lambda )|&\le  |\sin\varphi
|^{-4d-\frac72 n-6 }\ang\lambda ^{-2},\nonumber \\
|{\Cal K}_{R^{(1)}_MPQ_\lambda }(x,y,\lambda )|&\le  |\sin\varphi
|^{-4d-\frac72 n-7 }\ang\lambda ^{-2},\nonumber \\
|{\Cal K}_{R^{(2)}_MPQ_\lambda }(x,y,\lambda )|&\le  |\sin\varphi
|^{-\frac32 n-4 }\ang\lambda ^{-2}, \label{3.22a}
\end{align}

\end{thm}

\begin{proof}
 For the use of
Proposition \ref{pro5.1} we note that when  $T$ is
a $\psi $do of order $\le -n-1$, then 
\begin{align*}
\|T\|_{0,n+1}&\leg {\sum}_{|\gamma |\le n+1}\|D^\gamma T\|_{0,0},\\
\|T^*\|_{0,n+1}&=\|T\|_{-n-1,0}\eg \|T\ang D ^{n+1}\|_{0,0}.
\end{align*}
Consider $R^{(1)}_M$. The $N$ occurring in Proposition \ref{pro3.5}
can be written $N=\frac n2+\delta $, $\delta =\frac12$ or $1$. For $|\alpha |\le N$ and $|\beta |\le 1$, the symbols $D_x^\beta D_\xi ^\alpha (\xi ^\gamma
\circ r^{(1)}_M)$ are estimated by
$$
|D_x^\beta D_\xi ^\alpha (\xi ^\gamma \circ r^{(1)}_M)|\leg  |\sin\varphi |^{-2M-N-1-|\gamma |}\ang\xi ^{2d-M+|\gamma |}\ang{\lambda  }^{-2}.
$$
To apply Proposition \ref{pro5.1} with $|\gamma |$ up to $n+1$, we
must take the integer 
$M$ such that $2d-M+n+1\le 0$, so we let
$M=2d+n+1+\delta '$ with $\delta '\in [0,1[\,$. Then
$$
\|D^\gamma R^{(1)}_M\|_{0,0}\leg |\sin\varphi |^{-2M-1-N-|\gamma |}\ang{\lambda  }^{-2},
$$
and it follows that
\begin{align*}
\|R^{(1)}_M\|_{0,0}&\leg |\sin\varphi |^{-2M-1-N}\ang{\lambda
}^{-2}\\
\|R^{(1)}_M\|_{0,n+1}&\leg{\sum}_{|\gamma |\le n+1}\|D^\gamma
R^{(1)}_M\|_{0,0}\leg |\sin\varphi |^{-2M-N-n-2}\ang{\lambda
}^{-2}.
\end{align*}
Moreover,
$$
\|(R^{(1)}_M)^*\|_{0,n+1}=\|R^{(1)}_M\ang D^{n+1}\|_{0,0}\leg 
|\sin\varphi |^{-2M-1-N}\ang{\lambda
}^{-2}.$$
Insertion in \eqref{5.1} with $m=n+1$ gives that
\begin{align*}
|{\Cal K}_{R^{(1)}_M}&(x,y,\lambda )|\leg |\sin\varphi
|^{(-2M-N-n-2)\frac n{n+1}+(-2M-N-1)(1-\frac n {n+1})}\ang\lambda ^{-2}\\
&\leg |\sin\varphi
|^{-2M-N-n-1}\ang\lambda ^{-2}= |\sin\varphi
|^{-4d-\frac72 n-3-2\delta '-\delta }\ang\lambda ^{-2}\le  |\sin\varphi
|^{-4d-\frac72 n-6 }\ang\lambda ^{-2}.
\end{align*}
This shows the first estimate in \eqref{3.22a}.

For the second estimate we reuse the operator norms established for
$R^{(1)}_M$. They are now combined with some elementary norm estimates of
$PQ_\lambda $, namely:
$$
\|PQ_\lambda\|_{0,0}=\|1+\lambda Q_\lambda \|_{0,0}\leg 1+|\sin\varphi |^{-1} 
$$
holds  in view of \eqref{3.2b}, and moreover, for any $s\in{\Bbb R}$,
\begin{equation}
\|PQ_\lambda\|_{s,s}=\|P^{s/d}P Q_\lambda P^{-s/d}\|_{0,0} =\|P Q_\lambda \|_{0,0}\leg |\sin\varphi |^{-1} , \label{3.22b}
\end{equation}
since the operators commute. The composition with $PQ_\lambda $ thus
result in an extra factor $|\sin\varphi |^{-1}$ in the norm estimates,
hence likewise in the kernel estimate. This shows the second estimate
in \eqref{3.22a}.

For the third estimate we combine the elementary estimates
\begin{equation}
\|Q_\lambda \|_{s,s}\leg |\sin\varphi |^{-1}|\lambda |^{-1},\label{3.22c}
\end{equation}
shown earlier for $s=0$, and extendible to all $s$ by conjugation with
$P^{s/d}$,
with norm estimates of $R^{(2)}_MP$, derived similarly to above from the
last two lines in \eqref{3.21}. The latter give estimates in terms of
$|\sin\varphi |^{-1-N-1-n-1}\ang\lambda ^{-1}$. In the resulting
combined  estimate we can
replace $|\lambda |^{-1}$ by $\ang\lambda ^{-1}$, since we are working
under the hypothesis that a possible nullspace of $P$ has been removed.
This shows the last line in \eqref{3.22a}.
\end{proof}

Then we can show the estimate of the remainder kernel:

\begin{thm} \label{thm3.8}
Let $P$ be selfadjoint strongly elliptic of order
$d>0$ on $M$, with $\gamma (P)\ge 0$. The remainder kernel $\Cal
K_{V'_M}$ satisfies for $\arg
t=\theta $, $M=2d+ n+1+\delta '$ (with $\delta '\in [0,1[\,$): 
\begin{equation}
|\Cal K_{V'_M}(x,y,t)|\leg (\cos\theta 
)^{-2d-\frac 72 n-7}e^{-c'\operatorname{Re}t}|t|
,\label{3.27}
\end{equation}
where $c'>0$ if $\gamma (P)>0$,  $c'=0$ if $\gamma (P)=0$.
\end{thm}
\begin{proof} If $\gamma (P)>0$, we use the preceding estimates directly
to analyse
$$
\partial_t\Cal K_{V'_M}=-\tfrac i {2\pi }\int_{\Cal C}e^{-t\lambda
}\Cal K_{(Q _\lambda P)'_M}\,d\lambda ,\quad (Q_\lambda P)'_M =-R^{(1)}_M-R^{(2)}_M+R^{(1)}_MPQ_\lambda +R^{(2)}_MPQ_\lambda .
$$
The curve $\Cal C$ is chosen as in \eqref{3.1a}.
The first, third and fourth terms contribute with integrals of the form
$$
\tfrac i {2\pi }\int_{\Cal C}e^{-t\lambda
}f(x,y,\lambda )\,d\lambda 
$$
where 
$$
|e^{-\lambda t}f(x,y,\lambda )|\leg
e^{-c'\operatorname{Re}t}(\cos\theta )^{-4d-\frac 72 n-7}\ang\lambda ^{-2}
$$
on the curve by
Theorem \ref{thm3.5a} (recall that $\cos\theta \eg \sin\varphi _0$). 
Since $\ang\lambda ^{-2}$ integrates to $\infty
$, the resulting function is estimated by
$e^{-c'\operatorname{Re}t}(\cos\theta )^{-4d-\frac 72 n-7}$.

To find the contribution from  $R^{(2)}_M$, we first perform the
integration on the symbol level:
$$
\tfrac i {2\pi }\int_{\Cal C}e^{-t\lambda
}(p^0-\lambda )^{-1} p'_M\,d\lambda
=e^{-tp^0}p'_M.
$$
Here we can use an estimate from Lemma \ref{lem3.9}. By \eqref{3.33},
\begin{align*}
|D_x^\beta D_\xi ^\alpha (e^{-tp^0}p'_M)|&\leg (\sin \varphi_0)^{-M
-|\alpha|-|\beta|}
\ang\xi ^{-M-|\alpha|}e^{-c\operatorname{Re}t \ang \xi ^d}\\
&\leg (\sin \varphi_0)^{-M
-|\alpha|-|\beta|}
\ang\xi ^{-M-|\alpha|}e^{-c'\operatorname{Re}t }, \text{ hence}\\
|\Cal K_{\Op(e^{-p^0t}p'_M)}|&\leg(\sin \varphi_0)^{-M
}e^{-c'\operatorname{Re}t}\text{ when }M\ge n+1
.
\end{align*}

Since the latter estimate is  dominated
by that from the other terms, we conclude
that
$$
|\partial_t\Cal K_{V'_M}(x,y,t)|\leg
e^{-c'\operatorname{Re}t}(\cos\theta )^{-4d-\frac 72 n-7}.$$ 
 Then an
integration with respect to $t$ using that $\Cal K_{V'_M}(x,y,0)=0$
shows \eqref{3.27}.

In the case where $\gamma (P)=0$, the above considerations will be
valid for $V^\varepsilon (t)$ as in \eqref{3.19}. We then have to add
$(1-e^{-\varepsilon t})\Pi _0$, which has a smooth kernel bounded by
$\operatorname{min} \{|t|, 1\}$
and we reach the conclusion in the theorem.
\end{proof}

We can then finally show:

\begin{thm}\label{thm3.12}
Let $P$ be selfadjoint strongly elliptic of order
$d>0$ on $M$, with $\gamma (P)\ge 0$. The heat kernel $\Cal
K_{V}$ satisfies for all $t\in {\Bbb C}_+$ (with $\arg
t=\theta  \in \,]-\frac\pi 2,\frac\pi 2[\,$) the Poisson estimate,
where $N=\operatorname{max} \{\tfrac n d , \tfrac{7n}{2} + 4 d + 7
\}$:
\begin{equation}
|\Cal K_{V}(x,y,t)|\leg (\cos \theta  )^{-N}
e^{-\gamma(P)\operatorname{Re}t}\,\frac {|t|}{(d(x,y)+|t|^{1/d})^{d}}((d(x,y)+|t|^{1/d})^{-n}+1).\label{3.40}
\end{equation}
\end{thm}
\begin{proof}  
In local coordinates the estimate follows from Theorems \ref{thm3.11}
and \ref{thm3.8}, by choosing $M=n+1+2d+\delta '$ in Theorem \ref{thm3.8} and adding $\Cal K_{V_{-d-l}}$ for $0\leq l<M$; the most singular terms dominate.
This leads to the global estimate \eqref{3.40} (the effect of the
lower bound is handled as in Section 2).

\end{proof}

\begin{rem}\label{rem3.13}
Note that the order $d$ and dimension $n$ enter
 linearly in $N'=\tfrac{7n}{2} + 4 d + 7$,
and it is easy to see where the sizes come from: $2n+4d$ comes from the
power $-2M$ where $M\sim 2d+n$, one $n$ comes from the requirement for Agmon's estimate
Proposition \ref{pro5.1}, and $\frac n2 $ comes from the requirement
for Marschall's estimate Proposition \ref{pro3.5}. The number 7 includes
rounding up errors, and may be lowered. More substantial improvements would depend on choosing other
general principles; e.g.\ the $L_\infty $ kernel estimate in Beals
\cite{B70}, Lemma 2, which is slightly more efficient than Agmon's
estimate, and might save $\frac n2$ powers. Note that all orders
$d\in{\Bbb R}_+$ are allowed.

Our result applies in particular to the Dirichlet-to-Neumann operator
$P_{DN}$ of order $d=1$ associated with the Laplacian (and lower-order
perturbations of it).
For this operator, ter Elst and Ouhabaz \cite{EO13} have estimates in terms
of $-N''$-th powers of  $\cos\theta $, where the dimension $n$ enters
nonlinearly in $N$:
$$
N''=2n(n+1),\text{ compared to our }N'=\tfrac{7n}{2} + 11;
$$
here $N''>N'$ for $n\ge 6$.
They are proved by appealing to multiple
commutator estimates for 
semigroups defined from iterates of $P_{DN}$, refined $(L_p\to
L_q)$-estimates of Coifman-Meyer and others for pseudodifferential operators, Riesz potentials, interpolation, and
other tools. 
\end{rem}

\begin{rem} \label{rem3.16} Derivatives in $x$ and $y$ can also be
  estimated by these methods if needed. For example, $D_{x_j}\mathcal
  K_V(x,y,t)$ is described by the above formulas composed to the left
  with $D_{x_j}$; then in the remainder terms $D_{x_j}$ is composed
  with the $R^{(i)}_M$, giving symbols described in Proposition
  \ref{pro3.6}.
The remainder for $D_{y_j}\mathcal K_V(x,y,t)$ is described by
formulas, partly of the type where $R^{(i)}_M$ is composed to the
right with $D_{x_j}$ which just gives a factor $\xi _j$ on the symbol,
partly of a type containing $Q_\lambda D_{x_j} $. In the latter cases
one can use e.g.\ that
$$
PQ_\lambda D_{x_j}=PD_{x_j}Q_{\lambda }+P[D_{\xi _j},Q_\lambda ]
=D_{x_j}PQ_{\lambda }+[P,D_{x_j}]Q_{\lambda }
+PQ_\lambda [D_{x_j},P]Q_\lambda.
$$
Here the first term contributes a  $D_{x_j}$ that is absorbed in the $\psi
$do's before it, and in the other terms we note that $[D_{x_j},P]$ is
of order $d$ so that, by \eqref{3.22b},
\begin{align*}
\| [P,D_{x_j}]Q_\lambda \|_{0,0}&\leg \|PQ_\lambda \|_{0,0}\leg
|\sin\varphi |^{-1},\\
\|PQ_\lambda [D_{x_j},P]Q_\lambda \|_{0,0}&\leg \|PQ_\lambda \|_{0,0}\|[D_{x_j},P]Q_\lambda
\|_{0,0}\leg |\sin\varphi |^{-2}.
\end{align*}
\end{rem}

\section{Kernels of heat semigroups for
perturbations of
fractional Laplacians and the Dirichlet-to-Neumann operator}

This section complements the general upper bounds from Section 2
with lower estimates in the case of fractional powers of the Laplacian and the
Dirichlet-to-Neumann operator and their perturbations.

Let $\Delta $ be the (nonnegative) Laplace-Beltrami operator on the closed,
compact Riemannian $n$-dimensional manifold $M$.
The semigroups $e^{-t\Delta }$ and $V^d(t)= e^{-t\Delta^{d/2} \,}$ are related by subordination
formulas, which lead to an alternative proof of the upper kernel
estimate in this special case,  as well as to lower bounds. For $d=1$, they assume a simple form:

\begin{lem}\label{lem4.1}
Let $\lambda \ge 0$. One has for all $t\ge 0$:
\begin{equation}
 e^{-t\sqrt{\lambda }} = \frac{1}{2\sqrt{\pi }}\int_0^\infty e^{-s \lambda} t e^{-\frac{t^2}{4s}}
s^{-\frac32}ds\,.\label{4.1}
\end{equation}
\end{lem}
\begin{proof}  Let $\alpha =t\sqrt{\lambda }\,/2$ and let
$x=\frac{t}{2} s^{-\frac12}$; then $dx=-\frac{t}{4} s^{-\frac32}ds$, and
equation \eqref{4.1}
is turned into
\begin{equation}
\sqrt{\pi }\,e^{-2\alpha }=\int_0^\infty e^{-x^2-\frac{\alpha ^2}{x^2}}2\,dx
.\label{4.2}\end{equation}
To show this, note that the left-hand side $I(\alpha )$ satisfies
$
I(\alpha )\in C^1(\rp)$, $ \lim_{\alpha \to 0+}I(\alpha )=\sqrt{\pi }$,
and for $\alpha >0$ (with $y=\alpha x^{-1}$, $dy=-\alpha x^{-2}dx$):
$$
\partial_\alpha I(\alpha )=\int_0^\infty e^{-x^2-\frac{\alpha
^2}{x^2}}(-4\alpha  )x^{-2}\,dx=-2\int_0^\infty e^{-\frac{\alpha
^2}{y^2}-y^2}2\,dy
=-2I(\alpha ).
$$
Thus $I(\alpha )=ce^{-2\alpha }$ with $c=\sqrt{\pi }\,$.
\end{proof}

By Zolotarev \cite{Z86} (see also Grigor'yan \cite{G03}),  there
exists for any $0<d<2$ a
non-negative function
$\eta_t^d(s)$ such that
\begin{equation}
e^{-t\lambda^{d/2}} = \int_0^\infty e^{-s\lambda}\ \eta_t^d(s)\, ds  .\label{4.3}
\end{equation}
Here $\eta_t^d$ has the following properties
\begin{align}
\eta_t^d(s) &= t^{-2/d}\eta_1^d(\frac{s}{t^{2/d}}) \qquad (s,t>0)\ ,
\label{4.4}
\\
\eta_t^d(s) &\leg t s^{-1-\frac{d}{2}} \qquad (s,t>0)\ ,
\label{4.5}
\\
\eta_t^d(s) & \eg t s^{-1-\frac{d}{2}} \qquad (s\geq t^{2/d}>0) \
.\label{4.6}
\end{align}
By an application of the spectral theorem, we obtain for all $t>0$,
\begin{equation}
V^d(t)f = e^{-t\Delta^{d/2}} f = \int_0^\infty e^{-\tau\Delta}f \
\eta_t^d(\tau)\ d\tau \ , \text{ for
all }f\in H^s(M).\label{4.7}
\end{equation}
In view of \eqref{4.7}, it holds that
$$
\ang{\delta _x, V^d(t)\delta _y}=
\ang{\delta _x,\int_0^\infty e^{-\tau \Delta} \delta_y \, \eta_t^d(\tau)\, d\tau}
=\int_0^\infty \ang{\delta _x,e^{-\tau\Delta}\delta _y}\, \eta_t^d(\tau) \, d\tau,
$$
resulting in an identity for
the kernels: For all $t>0$,
\begin{equation}
{\Cal K}_{V^d \,}(x,y,t) = \int_0^\infty {\Cal K}_{e^{-\tau\Delta
}}(x,y)\ \eta_t^d(\tau) \,d\tau \ , \text{ for }(x,y)\in M\times M.\label{4.8}
\end{equation}
Using this formula, we can deduce upper and lower estimates for ${\Cal
K}_{V^d\,}$ from those known for ${\Cal K}_{e^{-\tau
{\Delta  }}}$.
The following upper and lower estimates are well-known (see e.g.\ L.\
Saloff-Coste \cite{S10}):
\begin{equation}
\frac{c_1}{{\Cal V}(x,\sqrt\tau \,)}e^{-C_1\frac{d(x,y)^2}\tau }
\le {\Cal K}_{e^{-\tau \Delta }}(x,y)\le
\frac{c_2}{\Cal{V}(x,\sqrt\tau \,)}e^{-C_2\frac{d(x,y)^2}\tau }.\label{4.9}
\end{equation}
Here $\Cal V(x,r)$ denotes the volume of a ball of radius $r$ around
$x$. For a closed compact $n$-dimensional manifold $M$, $\Cal V(x,r)\eg r^n$
for small $r$, and $\Cal V(x,r)$ equals the volume of the connected
component containing $x$ when $r\ge \operatorname{diam} M$. Hence
\begin{equation}
\Cal V(x,\sqrt\tau \,)^{-1}\eg (\tau ^{n/2})^{-1}+1.\label{4.10}
\end{equation}

\begin{thm}\label{thm4.2}
 Let $0<d<2$. The kernel of the semigroup
$V^d(t)=e^{-t\Delta^{d/2}}$ satisfies for all $t\ge 0$:
\begin{equation}
\Cal{K}_{e^{-t\Delta^{d/2}}}(x,y) \eg
\frac{t}{(d(x,y)+t^{1/d})^d}\left((d(x,y)+t^{1/d})^{-n}+1\right)\ . \label{4.11}
\end{equation}
\end{thm}

\begin{proof} The upper estimate follows already from Corollary
\ref{cor2.6}. The following proof moreover extends to give the lower estimate.
Inserting the
heat kernel bounds \eqref{4.9},
\eqref{4.10} into \eqref{4.8} and using \eqref{4.5}, we find
\begin{align*}
\Cal{K}_{V^d}(x,y,t) & \leg \int_0^\infty (\tau^{-n/2} + 1)\
\eta_t^d(\tau) \ e^{-C\frac{d(x,y)^2}{\tau}} d\tau \\
& \leg t \int_0^\infty \tau^{-n/2} \ \tau^{-1-\frac{d}{2}}\
e^{-C\frac{d(x,y)^2}{\tau}} d\tau + t \int_0^\infty  \tau^{-1-\frac{d}{2}}\ e^{-C
\frac{d(x,y)^2}{\tau}} d\tau\ .
\end{align*}
By a change of variables $\tau \mapsto C d(x,y)^2 \tau$, the first term equals
\begin{equation}t (Cd(x,y)^2)^{-\frac{d+n}{2}} \int_0^\infty \tau^{-\frac{n+d}{2}-1}\
e^{-1/\tau}\ d\tau \eg \frac{t}{d(x,y)^{n+d}}\ .\label{4.12}\end{equation}
Similarly, the second term is
$$t (Cd(x,y)^2)^{-\frac{d}{2}} \int_0^\infty \tau^{-\frac{d}{2}-1}\ e^{-1/\tau}\
d\tau \eg \frac{t}{d(x,y)^{d}}\ ,$$
and altogether,
$${\Cal K}_{V^d}(x,y,t)  \leg
\frac{t}{d(x,y)^{d}}\left(d(x,y)^{-n}+1\right) \ .$$
On the other hand, using the uniform bound ${\Cal K}_{e^{-\tau\Delta}}(x,y) \leg
\tau^{-n/2}+1$ and \eqref{4.4},
we obtain
\begin{align*}
{\Cal K}_{V^d}(x,y,t) & \leg \int_0^\infty (\tau^{-n/2} + 1)\
\eta_t^d(\tau)\ d\tau
 = \int_0^\infty (\tau^{-n/2}+1) \ \eta_1^d\left(\frac{\tau}{t^{2/d}}\right)\
t^{-2/d} \ d\tau \\
& = \int_0^\infty (t^{-n/d}\tau^{-n/2}+1) \ \eta_1^d(\tau)\ d\tau
 \eg t^{-n/d}+1 \ .
\end{align*}
Thus
$${\Cal K}_{V^d}(x,y,t) \leg \min\Bigl\{t^{-n/d}+1,
\frac{t}{d(x,y)^{d}}\left(d(x,y)^{-n}+1\right)\Bigr\} \ .$$
If $t^{1/d}\geq d(x,y)$,
$$t^{-n/d}\leg t^{-n/d}\Bigl(\frac{d(x,y)}{t^{1/d}}+1\Bigr)^{-n-d} =
t(d(x,y)+t^{1/d})^{-n-d}$$
and
$$1 \leg \Bigl(\frac{d(x,y)}{t^{1/d}}+1\Bigr)^{-d} =
t(d(x,y)+t^{1/d})^{-d} \ .$$
On the other hand, for $t^{1/d}\leq d(x,y)$ we have $d(x,y) \eg
d(x,y)+t^{1/d}$ and hence
$$\frac{t}{d(x,y)^{d}}\bigl(d(x,y)^{-n}+1\bigr)\leg
\frac{t}{(d(x,y)+t^{1/d})^{d}}\left((d(x,y)+t^{1/d})^{-n}+1\right) \ .$$
This shows ``$\leg$'' in \eqref{4.11}.

To show the opposite inequality in \eqref{4.11}, note that the integrand in \eqref{4.8}
is
non-negative, and \eqref{4.9}, \eqref{4.10} imply
$${\Cal K}_{V^d}(x,y,t) = \int_0^\infty
{\Cal K}_{e^{-\tau\Delta}}(x,y)\ \eta_t^d(\tau)\ d\tau \geg \int_\alpha^\infty (\tau^{-n/2} +
1)\ \eta_t^d(\tau) \ e^{-C\frac{d(x,y)^2}{\tau}} d\tau $$
for $\alpha = \max\{ t^{2/d}, d(x,y)^2\}$. Now, for $\tau \geq d(x,y)^2$,
$e^{-C\frac{d(x,y)^2}{\tau}} \geq e^{-C}$. Then by
\eqref{4.6},
\begin{align*}
&{\Cal K}_{V^d}(x,y,t)  \geg \int_\alpha^\infty (\tau^{-n/2} + 1)\ t
\tau^{1-\frac{1}{2}} \ d\tau
  \eg t \bigl(\alpha^{-\frac{n+d}{2}}+\alpha^{-\frac{d}{2}}\bigr) \\
&  = \min\bigl\{t^{-n/d}, t d(x,y)^{-n-d}\bigr\} + \min\bigl\{1, td(x,y)^{-d}\bigr\}
 \geq t (d(x,y)+t^{1/d})^{-n-d}+t(d(x,y)+t^{1/d})^{-d} \ .
\end{align*}
\end{proof}

For $d=1$ this complies well with the explicit kernel
formula \eqref{2.3} for the Poisson operator solving the Dirichlet problem
for the Laplacian on ${\Bbb R}^{n+1}_+$.

We also consider the case where $M$ is the boundary of  a compact
$(n+1)$-dimensional Riemannian manifold $\widetilde M$
with boundary. With $\Delta  $ denoting the nonnegative Laplace-Beltrami
operator on $M$, we shall compare  ${\Cal K}_{e^{-t\sqrt{\Delta  }\,}}$
with the kernel of the semigroup generated by the (nonnegative)
Dirichlet-to-Neumann
operator $P_{DN}$ on $M$. $P_{DN}$ is the operator mapping $u$ to the
normal derivative
$\partial_\nu \widetilde u$, where $\widetilde u$ is the harmonic
function on $\widetilde M$ with boundary value $u$.
It is known (cf.\ \cite{G71}) that $P_{DN}$ is an elliptic
pseudodifferential operator of order 1 on $M$ with the same principal
symbol as $\sqrt{\Delta  }\,$.

Since $\Delta ^{d/2}$ is a classical strongly elliptic $\psi $do of
order $d$,
Theorem \ref{thm2.5} applies to all operators of
the form $P=\Delta ^{d/2}+P'$ with $P'$ classical of order
$d-1$, giving upper estimates of the absolute value of the
kernels; note that no selfadjointness is required. For such operators
we can also show lower estimates.

\begin{thm} \label{thm4.3}
Let $d\in \,]0,2[\,$ and let $P$ be a classical $\psi
$do of order $d$ with the same principal
symbol as $\Delta^{d/2}$.
Then the kernel of $V(t)=e^{-tP}$ satisfies for all $t\ge 0$:
\begin{equation}
|\Cal K_{V}(x,y,t)|\leg  \frac{t}{(d(x,y)+t^{1/d})^{d}}\,\Bigl( \frac 1{(d(x,y)+t^{1/d})^{-n}}+1\Bigr)
+e^{-c_1t}\frac{t}{(d(x,y)+t^{1/d})^{d+n-1}},\label{4.13}
\end{equation}
for any $c_1<\gamma (P)$ ($c_1=\gamma (P)$ if Corollary {\rm 2.6}
applies). Moreover, there is an $r>0$ such that
\begin{equation}
|\Cal K_{V}(x,y,t)|\geg
 t\,(d(x,y)+t^{1/d})^{-d-n} ,\text{ for }d(x,y)+t^{1/d}\le r.\label{4.14}
\end{equation}
\end{thm}
\begin{proof} As $P$ and $\Delta^{d/2}$ have the same principal symbol,
$$
V(t)=V^d(t)+V' ,
$$
where $V'$ is of lower order, more precisely $V'$
is the difference between the first remainders for $V(t)=e^{-tP}$ and
$V^d(t)=e^{-t\Delta^{d/2}}$, as in the second line of \eqref{2.25}. Hence
\begin{equation}
|\Cal K_{V'}(x,y,t)|\leg  e^{-c_1 t}t\, (d(x,y)+t^{1/d})^{1-n-d}.
\label{4.15}
\end{equation}
Now \eqref{4.11} and \eqref{4.15} together imply \eqref{4.13}.

To obtain the lower estimate \eqref{4.14}, we note that \begin{equation}
cs^{-n-d}-c's^{1-n-d}=cs^{-n-d}(1-c'c^{-1}s)\ge 2^{-1}cs^{-n-d},\text{
when }s\le c/(2c'),\label{4.16}
\end{equation}
so for all $t$ in a bounded
set where $e^{-c_1t}\le c'$, the lower estimate in
\eqref{4.11} implies that \eqref{4.14} holds for all small $d(x,y)+t^{1/d}$.
\end{proof}

We can also obtain upper and lower estimates for the Dirichlet-to-Neumann
 operator.

\begin{thm}\label{thm4.4}
The kernel of $e^{-tP_{DN}}$ satisfies for all $t\ge 0$:
\begin{equation}
\Cal K_{e^{-tP_{DN}}}(x,y,t)\leg
\frac{t}{d(x,y)+t}\,\left((d(x,y)+t)^{-n}+1\right),\label{4.17}
\end{equation}
and there is an $r>0$ such that it satisfies
\begin{equation}
\Cal K_{e^{-tP_{DN}}}(x,y,t)\geg t\,(d(x,y)+t)^{-1-n},\text{ for
}d(x,y)+t\le r.\label{4.18}
\end{equation}
\end{thm}

\begin{proof} Here $P_{DN}$ is known to be selfadjoint nonnegative, and
the semigroup has real,
nonnegative kernel (\cite{AM07, AM12}), so that we may omit absolute values.
The upper estimate \eqref{4.17} follows from Corollary 2.6. The lower
estimate \eqref{4.18} follows from Theorem \ref{thm4.3} since $P_{DN}$ differs from $\Delta
^{1/2}$ by a
classical $\psi $do of order 0.
\end{proof}

\begin{rem} This work was inspired from a conversation of the second
  author with W.\ Arendt
  and A.\ ter Elst in August 2012 on the need for kernel estimates for
  the Dirichlet-to-Neumann semigroup, where
we suggested the applicability of pseudodifferential methods as in
\cite{G96}. When our first version of this paper was posted in
arXiv:1302.6529,
we learned of the efforts of ter Elst and Ouhabaz
in \cite{EO13}, giving an analysis of the Dirichlet-to-Neumann
semigroup kernel by
somewhat different methods, and obtaining some of the same results as
those presented here.

\end{rem}

\end{document}